\documentclass[preprint,3p]{elsarticle}
\usepackage[toc,page,title,titletoc,header]{appendix}
\usepackage{stmaryrd}
\SetSymbolFont{stmry}{bold}{U}{stmry}{m}{n}
\usepackage{amsmath}
\usepackage{amssymb}
\usepackage{algorithm}
\usepackage{amsthm}
\usepackage{tabularx}
\usepackage{subfigure}
\usepackage{epsfig}
\usepackage{indentfirst}
\usepackage{cases}
\biboptions{sort&compress}
\usepackage{graphicx}
\usepackage{epstopdf}
\usepackage{color}
\usepackage{float}
\usepackage{algorithmic}
\usepackage{bm}
\usepackage{subfig}
\usepackage{multirow}

\newcommand{\be}{\begin{equation}}
\newcommand{\ee}{\end{equation}}
\newcommand{\ba}{\begin{array}}
\newcommand{\ea}{\end{array}}
\newcommand{\bea}{\begin{eqnarray}}
\newcommand{\eea}{\end{eqnarray}}
\newcommand{\beas}{\begin{eqnarray*}}
\newcommand{\eeas}{\end{eqnarray*}}

\newtheorem{thm}{Theorem}[section]

\newtheorem{remark}{Remark}[section]

\newtheorem{ex}{Example}[section]

\numberwithin{equation}{section}

\def\R{{\mathbb R}}
\def\S{{\mathbb S}}

\newcommand{\p}{\partial}

\newcommand{\cE}{\mathcal E}

\newcommand{\eps}{\varepsilon}

\renewcommand{\l}{\left}
\renewcommand{\r}{\right}

\renewcommand{\d}{\mathrm{d}}

\newcommand{\bn}{\mathbf{n}}

\newcommand\cV{{\mathcal V}}
\newcommand{\bX}{\mathbf{X}}

\newcommand{\bZ}{\mathbf{Z}}

\begin{document}

\begin{frontmatter}

\title{Predictor-corrector, BGN-based parametric finite element methods for surface diffusion}

\author[1,4]{Wei Jiang}
\address[1]{School of Mathematics and Statistics, Wuhan University, Wuhan, 430072, China}
\address[4]{Shenzhen Research Institute of Wuhan University, Shenzhen, 518057, China}
\ead{jiangwei1007@whu.edu.cn}

\author[2]{Chunmei Su}
\ead{sucm@tsinghua.edu.cn}
\author[2]{Ganghui Zhang}
\ead{gh-zhang19@mails.tsinghua.edu.cn}
\address[2]{Yau Mathematical Sciences Center, Tsinghua University, Beijing, 100084, China}
\author[1]{Lian Zhang}
\ead{lian\_zhang@whu.edu.cn}


\begin{abstract}

We present a novel parametric finite element approach for simulating the surface diffusion of curves and surfaces. Our core strategy incorporates a predictor-corrector time-stepping method, which enhances the classical first-order temporal accuracy to achieve second-order accuracy. Notably, our new method eliminates the necessity for mesh regularization techniques, setting it apart from previously proposed second-order schemes by the authors (J. Comput. Phys. 514 (2024) 113220). Moreover, it maintains the long-term mesh equidistribution property of the first-order scheme. The proposed techniques are readily adaptable to other geometric flows, such as (area-preserving) curve shortening flow and surface diffusion with anisotropic surface energy. Comprehensive numerical experiments have been conducted to validate the accuracy and efficiency of our proposed methods, demonstrating their superiority over previous schemes.
\end{abstract}




\end{frontmatter}
\section{Introduction}

Surface diffusion plays an important role in the field of materials science. It was initially proposed by Mullins in 1957 as a mathematical model to describe the evolution of grain boundaries in polycrystal materials~\cite{Mullins1957}. Since then,
this pioneering work has sparked extensive research and application across multiple research areas, such as the crystal growth \cite{Gilmer1972, Gomer1990}, the evolution of voids in microelectronic circuits \cite{Li1999},
solid-state dewetting \cite{Wang2015, Jiang2020, Zhao2020, Li2023}, and computational geometry for surface design \cite{Schneider2001}.

This paper focuses on the numerical simulation of surface diffusion flow (SDF) and its extension to other types of geometric flows. We begin with the planar surface diffusion flow, investigating the evolution of a family of simple closed curves denoted as $\Gamma(t)$. The evolution is governed by the following geometric equation:
  \begin{equation}\label{SDF:geo equation}
	 \cV = \p_{ss}\kappa\ \bn,
  \end{equation}
where $\cV$ represents the velocity, $s:=s(t)$ is the arc-length, $\kappa$ denotes the plane curvature of the curve, and $\bn$ is the outward unit normal vector.  It is widely known that the above SDF can be interpreted as the $H^{-1}$-gradient flow of the perimeter functional \cite{Barrett2020,Garcke2013}.

The simulation of the SDF and its anisotropic variant has attracted extensive research interest, leading to the development of numerous numerical methods. Notable examples include the graph formulation \cite{Baensch2004,Xu2009}, the finite difference method \cite{Mayer2001}, the $\theta$-$L$ formulation \cite{Huang2021,Huang2024}, and the finite element method based on a complicated variational formulation \cite{Baensch2005}.
In their pioneering work \cite{Barrett2007a}, Barrett, Garcke and N\"urnberg (referred to as BGN) introduced an unconditionally energy-stable parametric finite element method (PFEM) for isotropic SDF by parametrizing $\Gamma(t)$ using $\bX(\cdot,t):\mathbb{I}=[0,1]\rightarrow \R^2$ and introducing a novel coupled reformulation of the SDF \eqref{SDF:geo equation} as follows:
\begin{equation}\label{SDF:Coupled equation}
\begin{split}
	\p_t \mathbf{X}\cdot \mathbf{n} &= \p_{ss}\kappa,\\
  	\kappa\, \mathbf{n}&=-\p_{ss}\mathbf{X}.
\end{split}
\end{equation}
Noteworthy, the PFEM based on the above BGN formulation exhibits excellent mesh quality, characterized by asymptotic long-term mesh equidistribution~\cite{Barrett2007a, Barrett2020}. They further extended this approach to the anisotropic SDF with Riemannian-like anisotropic surface energy $\gamma(\bn)$  \cite{Barrett2008a}. Similar formulation has been applied to other geometric flows, such as curve shortening flow (CSF) and area-preserving curve shortening flow (AP-CSF) \cite{Barrett2007}. More recently, Bao, Jiang and Li utilized the  Cahn-Hoffman $\bm{\xi}$-vector formulation \cite{Jiang2019} to address anisotropic surface diffusion flow (A-SDF) with arbitrary anisotropic surface energy densities. However, it is important to note that these numerical methods utilize the backward Euler method for time discretization, resulting in only first-order temporal accuracy. For numerical methods with higher-order temporal accuracy in simulating isotropic SDF, we refer to the recent studies \cite{Hu2022,Duan2023,Kovacs2021,Jiang2023,Jiang2024}.

In this paper, we present a novel strategy to enhance the temporal accuracy of the previously discussed first-order PFEMs. This approach is applicable to the evolution of curves and surfaces driven by isotropic SDF, CSF/AP-CSF, and A-SDF. Our methodology incorporates a predictor-corrector technique for discretizing the time variable, thereby improving the temporal accuracy.
Specifically, when computing  the numerical solution from time $t$ to time $t+\tau$, we first utilize a first-order scheme for a half time step to get an approximation of the solution at time $t+\frac{\tau}{2}$. Subsequently, based on this `predictor', we employ a formally second-order scheme with a local truncation error of $\mathcal{O}(\tau^2)$ to derive an approximation of the solution at time $t+\tau$.

%
\smallskip
The predictor-corrector strategy introduced in this paper offers several notable advantages compared to previous methodologies:
\begin{itemize}
\item {\textbf{Second-Order Accuracy}}: It demonstrates robust second-order accuracy, as assessed using the shape metric referenced in \cite{Jiang2023}. This represents a significant improvement in temporal accuracy compared to previous methods.
\item {\textbf{Retained Favorable Mesh Properties}}: The strategy retains the beneficial mesh properties of the traditional first-order BGN scheme, which eliminates the need for mesh regularization. This contrasts with the authors' earlier work \cite{Jiang2023}, which employed the Crank-Nicolson leap-frog discretization.

\item {\textbf{Versatility and Extensibility}}: The method can be readily adapted to various types of surface diffusion and, more broadly, to other geometric flows. Crucially, this extension maintains the same second-order accuracy and advantageous mesh quality, ensuring consistency and reliability across diverse applications.
\end{itemize}

Throughout the paper, we will primarily focus on the applications of several key geometric flows in addition to the SDF. These include, but are not limited to the following:
\begin{itemize}
	\item The CSF and AP-CSF, for which the motion is governed by the equation
	\begin{equation}\label{CSF/AP-CSF:geo equation}
	 \cV =\begin{cases}
	 	-\kappa  \bn, \quad\qquad& \text{CSF},\\
	 	-\l(\kappa-\l<\kappa\r>\r)\bn, \quad\qquad & \text{AP-CSF}.
	 \end{cases}
  \end{equation}
  Similar to \eqref{SDF:Coupled equation}, the corresponding BGN formulation \cite{Barrett2007,Barrett2020} reads as:
  \begin{equation}\label{CSF/AP-CSF:Coupled equation}
\begin{split}
	\p_t \mathbf{X}\cdot \mathbf{n} &= \begin{cases}
		-\kappa,\qquad& \text{CSF},\\
	   -\l(\kappa-\l<\kappa\r>\r),\qquad & \text{AP-CSF},
	\end{cases}\\
  	\kappa\, \mathbf{n}&=-\p_{ss}\mathbf{X}.
\end{split}
\end{equation}
	\item The A-SDF, for which the motion is governed by the equation:
	\begin{equation}\label{aSDF:geo equation}
	 \cV = \p_{ss}\mu\ \bn,
  \end{equation}
  where $\mu$ is the chemical potential \cite{Cahn94,Bao17}. From a variational point of view, it can be interpreted as an $H^{-1}$-gradient flow of the energy  functional $W:=\int_{\Gamma}\gamma(\bn)\ \d s$, where $\gamma(\bn):\S^1\rightarrow \R^+$ denotes the anisotropic surface energy density \cite{Bao2022a,Bao17,Barrett2020}. Following the work of Bao, Jiang and Li  \cite{Bao2022a,Bao2022b}, a coupled reformulation can be expressed as
\begin{equation}\label{aSDF:Coupled equation1}
\begin{split}
	\p_t \mathbf{X}\cdot \mathbf{n} &= \p_{ss}\mu,\\
  	\mu\, \mathbf{n}&=-\p_{s}\bm{\xi}^{\perp},
\end{split}
\end{equation}
where $\bm{\xi}=\l.\nabla\gamma(\mathbf{p})\r|_{\mathbf{p}=\bn}$ represents the Cahn-Hoffman $\bm{\xi}$-vector \cite{Cahn98, Jiang2019}, and the notation $\cdot^\perp$ indicates a clockwise rotation by $\frac{\pi}{2}$.
Subsequently, by introducing a stabilizing function $k(\bn):\S^1\rightarrow\R^+$ which depends on the anisotropic surface energy $\gamma(\bn)$, and a positive symmetrized matrix $\bZ_{k}(\bn)$ depending on $k(\bn)$, \eqref{aSDF:Coupled equation1} could be further reformulated as
\begin{equation}\label{aSDF:Coupled equation2}
\begin{split}
	\p_t \mathbf{X}\cdot \mathbf{n} &= \p_{ss}\mu,\\
  	\mu\, \mathbf{n}&=-\p_{s}(\bZ_{k}(\bn)\p_s\bX).
\end{split}
\end{equation}
	\item The SDF of surfaces in $\mathbb{R}^3$, for which the motion is governed by the equation
  \begin{equation}\label{SDF3d:geo equation}
	 \cV = \Delta_{\Gamma}H \ \bn,
  \end{equation}
  where $H$ is the mean curvature of the surface, and $\Delta_{\Gamma}$ is the surface Laplacian operator on the surface $\Gamma$. The corrsponding BGN formulation is given by
  \begin{equation}\label{SDF3d:Coupled equation}
\begin{split}
	\p_t \mathbf{X}\cdot \mathbf{n} &= \Delta_{\Gamma} H,\\
   H\, \mathbf{n}&= -\Delta_{\Gamma}\mathrm{Id},
\end{split}
\end{equation}
where $\mathrm{Id}$ is the identity operator defined on $\Gamma$.
\end{itemize}

The remainder of the paper is organized as follows. Section 2 begins with a review of the classical first-order BGN scheme for solving planar isotropic SDF and then introduces a semi-implicit temporal second-order scheme to enhance the accuracy. We also establish the well-posedness of this scheme and prove its long-time mesh equidistribution property. In Section 3, we extend this methodology to CSF, AP-CSF, A-SDF and SDF in $\R^3$. Section 4 includes extensive numerical experiments that demonstrate the accuracy and efficiency of our proposed methods across various types of geometric flows.
Finally, we draw some conclusions and discuss potential directions for future research in Section 5.

\section{A predictor-corrector second-order algorithm for planar isotropic surface diffusion flow}

We first recall the classical first-order BGN scheme for planar isotropic SDF \cite{Barrett2007a,Barrett2020}. As detailed in the introduction, the isotropic  SDF is reformulated as \eqref{SDF:Coupled equation}. For finite element approximation, let $\mathbb{I}=[0,1]= \bigcup_{j=1}^N I_j$, $N\ge 3$, be a partition of $\mathbb{I}$ with intervals $I_j=[\rho_{j-1},\rho_j]$. We define the linear finite element space as
\[
V^h:=\{u\in C(\mathbb{I}): u|_{I_j} \,\,\, \mathrm{is\,\,\,linear,\,\,\,} \forall j=1,2,\ldots,N;\quad u(\rho_0)=u(\rho_N) \}\subseteq H^1(\mathbb{I}).
\]
The mass lumped inner product $(\cdot,\cdot)_{\Gamma^h}^h$ over the polygonal curve $\Gamma^h:=\bX^h\in [V^h]^2$, which is an approximation of $(\cdot,\cdot)_{\Gamma^h}$ by using the composite trapezoidal rule, is defined as
\[
(u,v)_{\Gamma^h}^h:=\frac{1}{2}\sum_{j=1}^N|\bX^h(\rho_j,t)-\bX^h(\rho_{j-1},t)|\l[(u\cdot v)(\rho_j^-)+(u\cdot v)(\rho_{j-1}^+) \r],
\]
where $u, v$ are two scalar/vector piecewise continuous functions with possible jumps at the nodes $\{\rho_j\}_{j=1}^N$,
and $u(\rho_j^{\pm})=\lim\limits_{\rho\rightarrow \rho_j^{\pm}}u(\rho)$.

Subsequently, the semi-discrete scheme for the formulation \eqref{SDF:Coupled equation} is as follows: given the initial polygon $\Gamma^h(0)$ with vertices lying on the initial curve $\Gamma(0)$ in a clockwise orientation, parametrized by $\bX^h(\cdot,0)\in [V^h]^2$,
find $(\bX^h(\cdot,t),\kappa^h(\cdot,t))\in [V^h]^2\times V^h$ such that
\begin{equation}\label{SDF:Semi-discrete}
	\begin{split}
			&\l(\p_t\mathbf{X}^h\cdot  \mathbf{n}^h,\varphi^h \r)_{\Gamma^h}^h+\l(  \p_s\kappa^h,\p_s\varphi^h\r)_{\Gamma^h}=0,\quad \forall\ \varphi^h\in V^h,\\
			&\l(\kappa^h,\mathbf{n}^h\cdot \bm{\omega}^h\r)^h_{\Gamma^h}-\l(\p_s \mathbf{X}^h,\p_s\bm{\omega}^h \r)_{\Gamma^h}=0,\quad \forall\ \bm{\omega}^h\in  [V^h]^2,
		\end{split}
\end{equation}
where we always integrate over the current curve $\Gamma^h$ described by $\mathbf{X}^h$, the outward unit normal $\bn^h$ is a piecewise constant vector given by
\[\bn^h|_{I_j}=-\mathbf{h}_j^\perp/|\mathbf{h}_j|, \quad \mathbf{h}_j=\bX^h(\rho_j,t)-\bX^h(\rho_{j-1},t),\quad j=1,\ldots, N,\]
 and the partial derivative $\p_s$ is defined piecewisely over each side of the polygon
$$\p_s f|_{I_j}=\frac{\p_\rho f}{|\p_\rho \mathbf{X}^h|}|_{I_j}=\frac{(\rho_j-\rho_{j-1})\p_\rho f|_{I_j}}{|\mathbf{h}_j|}.$$ It was shown that the scheme \eqref{SDF:Semi-discrete} always equidistributes the vertices along $\Gamma^h$ for $t>0$ if they are not locally parallel \cite{Barrett2007a}. 

For the fully discrete scheme, we select a uniform time step size $\tau>0$ for simplicity. Let $\bX^m\in [V^h]^2$ and $\Gamma^m$ be the approximations of $\bX(\cdot,t_m)$ and $\Gamma(t_m)$, respectively, for $m=0,1,2,\ldots$, with $t_m:=m\tau$. Define $\mathbf{h}_j^m:=\bX^m(\rho_j)-\bX^m(\rho_{j-1})$ and assume that $|\mathbf{h}_j^m|>0$ for $j=1,\ldots,N$ and $m>0$. The discrete unit normal vector $\bn^m$, the discrete inner product $(\cdot,\cdot)^h_{\Gamma^m}$ and the discrete operator $\p_s$ are defined analogously to their semi-discrete counterparts. In view of the formal first-order approximation for $\p_t \bX$, $\kappa$ and $\p_s\bX$:
\begin{align*}
	&\p_t \bX(\cdot, t_m)
	= \frac{\mathbf{X}(\cdot,t_{m+1})-\mathbf{X}(\cdot, t_m)}{\tau}+\mathcal{O}(\tau),\\
	&\p_s\kappa(\cdot,t_m)
	=\p_s\kappa(\cdot,t_{m+1})+\mathcal{O}(\tau),\quad
	\p_s\bX(\cdot,t_m)= \p_s \mathbf{X}(\cdot, t_{m+1})+\mathcal{O}(\tau),
\end{align*}
Barrett, Garcke and N\"urnberg \cite{Barrett2007,Barrett2007a} proposed the following semi-implicit  scheme (referred to as the BGN scheme).

\noindent
(\textbf{First-order BGN scheme}):~For $m\ge 0$, find $\mathbf{X}^{m+1}\in [V^h]^2$ and $\kappa^{m+1}\in V^h$ such that
\begin{equation}\label{SDF:BGN1}
		\begin{split}
			&\l(\frac{\mathbf{X}^{m+1}-\mathbf{X}^m}{\tau},\varphi^h \mathbf{n}^m \r)^h_{\Gamma^m}+\l( \p_s \kappa^{m+1},\p_s \varphi^h \r)_{\Gamma^m}=0,\quad \forall\ \varphi^h\in V^h,\\	
			&\l(\kappa^{m+1},\mathbf{n}^m\cdot \bm{\omega}^h\r)_{\Gamma^m}^h-\l(\p_s \mathbf{X}^{m+1},\p_s\bm{\omega}^h\r)_{\Gamma^m}=0,\quad \forall\ \bm{\omega}^h\in  [V^h]^2.
		\end{split}
	\end{equation}
The well-posedness and unconditional energy stability have been established under mild conditions. In practice, the BGN scheme \eqref{SDF:BGN1} demonstrates quadratic convergence in space \cite{Barrett2007a} and linear convergence in time \cite{Jiang2023}. Furthermore, as discussed in \cite[Section 4.6]{Barrett2020}, the BGN scheme exhibits favorable mesh quality and an asymptotic long-time mesh equidistribution property in practical applications.

\smallskip

 We now clarify our predictor-corrector strategy for the planar isotropic SDF and introduce a second-order temporal scheme based on this approach and the formulation \eqref{SDF:Coupled equation}. We establish the fundamental properties of the scheme, including well-posedness and the long-time equidistribution property. Specifically, we approximate the  terms $\p_t\bX$, $\kappa$ and $\p_s \bX$ using formal Taylor expansion formulas as follows:
\begin{equation}\label{PC method approximation}
\begin{split}
	\p_t \bX(\cdot,t_{m+\frac{1}{2}}) &= \frac{\mathbf{X}(\cdot,t_{m+1})-\mathbf{X}(\cdot,t_m)}{\tau}+\mathcal{O}(\tau^2),\\
	\p_{s}\kappa(\cdot,t_{m+\frac{1}{2}})
	&= \frac{\p_s\kappa(\cdot,t_{m+1})+\p_s\kappa(\cdot,t_{m})}{2}+\mathcal{O}(\tau^2),\\
	\p_s\bX(\cdot,t_{m+\frac{1}{2}})&= \frac{\p_s \mathbf{X}(\cdot, t_{m+1})+\p_s \mathbf{X}(\cdot, t_{m})}{2}+\mathcal{O}(\tau^2).
\end{split}	
\end{equation}
Inspired by the above approximations, we derive a formal second-order scheme as follows:
\begin{equation}\label{SDF:BGNCN}
		\begin{split}
	&\l(\frac{\mathbf{X}^{m+1}-\mathbf{X}^{m}}{\tau},\varphi^h \mathbf{n}^{m+\frac{1}{2}} \r)^h_{\Gamma^{m+\frac{1}{2}}}+\l( \p_s\l(\frac{\kappa^{m+1}+\kappa^{m}}{2}\r),\p_s\varphi^h \r)_{\Gamma^{m+\frac{1}{2}}}=0,\quad \forall\ \varphi^h\in V^h,\\
&\l(\frac{\kappa^{m+1}+\kappa^{m}}{2},\mathbf{n}^{m+\frac{1}{2}}\cdot \bm{\omega}^h\r)_{\Gamma^{m+\frac{1}{2}}}^h-\l(\p_s\l(\frac{ \mathbf{X}^{m+1}+ \mathbf{X}^{m}}{2}\r),\p_s\bm{\omega}^h\r)_{\Gamma^{m+\frac{1}{2}}}=0, \quad \forall\ \bm{\omega}^h\in  [V^h]^2,
		\end{split}
	\end{equation}
where $\Gamma^{m+\frac{1}{2}}$ is a suitable approximation of $\Gamma(t_{m+1/2})$, e.g., $\Gamma^{m+\frac{1}{2}}:=\bX^{m+\frac{1}{2}}=\frac{\bX^m+\bX^{m+1}}{2}$. To obtain a semi-implicit scheme that is easier to solve, we compute the `predictor' $\widetilde{\bX}^{m+\frac{1}{2}}$ as the solution of the classical BGN scheme \eqref{SDF:BGN1}, using the initial data $\bX^m$ and time step of $\tau/2$. We denote $\widetilde{\Gamma}^{m+\frac{1}{2}}:=\widetilde{\bX}^{m+\frac{1}{2}}$. Let $\widetilde{\bn}^{m+\frac{1}{2}}$ be the normal vector of the polygon $\widetilde{\bX}^{m+\frac{1}{2}}$. By substituting this predictor into \eqref{SDF:BGNCN}, we obtain the following BGN/PC scheme.

\smallskip

\noindent(\textbf{Second-order BGN/PC scheme}):~ For $\bX^0 \in [V^h]^2$ and $\kappa^0\in  V^h$, which serve as suitable approximations at the initial time level $t_0=0$, for $m\ge 0$, let $\widetilde{\Gamma}^{m+\frac{1}{2}}:=\widetilde{\bX}^{m+\frac{1}{2}}\in [V^h]^2$ be the solution of the classical BGN scheme \eqref{SDF:BGN1} with initial data $\bX^m$ and a time step of $\tau/2$. The objective is to find $\mathbf{X}^{m+1}\in [V^h]^2$ and $\kappa^{m+1}\in V^h$ such that
\begin{equation}\label{SDF:BGN/PC}
		\begin{split}
	&\l(\frac{\mathbf{X}^{m+1}-\mathbf{X}^{m}}{\tau},\varphi^h \widetilde{\mathbf{n}}^{m+\frac{1}{2}} \r)^h_{\widetilde{\Gamma}^{m+\frac{1}{2}}}+\l( \p_s\l(\frac{\kappa^{m+1}+\kappa^{m}}{2}\r),\p_s\varphi^h \r)_{\widetilde{\Gamma}^{m+\frac{1}{2}}}=0,\quad \forall\ \varphi^h\in V^h,\\	&\l(\frac{\kappa^{m+1}+\kappa^{m}}{2},\widetilde{\mathbf{n}}^{m+\frac{1}{2}}\cdot \bm{\omega}^h\r)_{\widetilde{\Gamma}^{m+\frac{1}{2}}}^h-\l(\p_s\l(\frac{ \mathbf{X}^{m+1}+ \mathbf{X}^{m}}{2}\r),\p_s\bm{\omega}^h\r)_{\widetilde{\Gamma}^{m+\frac{1}{2}}}=0, \quad \forall\ \bm{\omega}^h\in  [V^h]^2.
		\end{split}
	\end{equation}
Here, the normal vector is defined as $\widetilde{\mathbf{n}}^{m+\frac{1}{2}}:=-\l(\frac{\p_\rho \widetilde{\bX}^{m+\frac{1}{2}}  }{|\p_\rho \widetilde{\bX}^{m+\frac{1}{2}} |}\r)^{\perp}$, and the derivative $\p_s$ refers to the arc length of $\widetilde{\Gamma}^{m+\frac{1}{2}}$. We present the following BGN/PC algorithm, which combines the schemes \eqref{SDF:BGN1} and \eqref{SDF:BGN/PC}.

\begin{algorithm}\textbf{BGN/PC algorithm for planar isotropic case}
\begin{algorithmic}[1]\label{BGN/PC algorithm}
 \REQUIRE An initial curve $\Gamma(0)$ approximated by a polygon $\Gamma^0$ with $N$ vertices, described by $\bX^0\in [V^h]^2$, time step $\tau$ and terminate time $T$ satisfying $T/\tau\in\mathbb{N}$.
 \ENSURE Computational solution $\Gamma^{T/\tau}:=\bX^{T/\tau}$.
 \STATE Calculate $\kappa^0$ by  $\Gamma^0$ and a formula based on formula $\kappa\bn=-\p_{ss}\bX$ and least square method \cite{Barrett2007a,Jiang2023}.   Set $m=0$.
 \WHILE{$m < T/\tau$,}
\STATE {\small
 \begin{itemize}
\item \textbf{Predictor}: Compute  $\widetilde{\Gamma}^{m+\frac{1}{2}}:= \widetilde{\bX}^{m+\frac{1}{2}}$ by using BGN scheme \eqref{SDF:BGN1} with $\Gamma^m$ and $\tau/2$.
 \item \textbf{Corrector}: Compute $\Gamma^{m+1}:=\bX^{m+1},\kappa^{m+1}$ by using the BGN/PC scheme \eqref{SDF:BGN/PC} with
 $\bX^{m},\kappa^{m},\widetilde{\Gamma}^{m+\frac{1}{2}}$ and $\tau$.
\end{itemize}
}
 \STATE $m = m + 1$;
 \ENDWHILE
 \end{algorithmic}
\end{algorithm}

\smallskip

\begin{remark}\label{Remark:BGN2}
 In our previous work \cite{Jiang2023}, we introduced a Crank-Nicolson leap-frog (CNLF) time-stepping discretization and proposed the following scheme, referred to here as the BGN/CNLF scheme:

\noindent(\textbf{Second-order BGN/CNLF scheme}):~ For $(\bX^0,\kappa^0), (\bX^1,\kappa^1)\in [V^h]^2\times V^h$, which represent appropriate approximations at the time levels $t_0=0$ and $t_1=\tau$, find $(\mathbf{X}^{m+1}, \kappa^{m+1})\in [V^h]^2\times V^h$ for $m\ge 1$ such that
\begin{equation}\label{SDF:BGN2}
\begin{split}
&\l(\frac{\mathbf{X}^{m+1}-\mathbf{X}^{m-1}}{2\tau},\varphi^h \mathbf{n}^m \r)^h_{\Gamma^m}+\l(  \frac{\p_s\kappa^{m+1}+\p_s\kappa^{m-1}}{2},\p_s\varphi^h \r)_{\Gamma^m}=0,\quad \forall\ \varphi^h\in V^h,\\
&\l(\frac{\kappa^{m+1}+\kappa^{m-1}}{2},\mathbf{n}^m\cdot \bm{\omega}^h\r)_{\Gamma^m}^h-\l(\frac{\p_s \mathbf{X}^{m+1}+\p_s \mathbf{X}^{m-1}}{2},\p_s\bm{\omega}^h\r)_{\Gamma^m}=0,\quad \forall\ \bm{\omega}^h\in  [V^h]^2.
\end{split}
\end{equation}
The derived two-step scheme is well-posed and demonstrates second-order accuracy in time. However, mesh regularization may be necessary in certain scenarios to prevent distortion, due to some oscillatory behavior. The proposed BGN/PC algorithm  effectively addresses the oscillatory behavior of the mesh ratio (see Example \ref{Evolution of geometric quantities} and Figure \ref{Fig:GEO}). This is attributed to the incorporation of a regularly updated mesh at each time step, utilizing the classical BGN scheme, thereby maintaining desirable mesh properties. Additionally, the order of the local truncation error ensures consistent second-order accuracy.
\end{remark}

\smallskip

\begin{remark}
It would also be beneficial to compare the BGN/PC scheme with another second-order BGN-based scheme introduced in the authors' previous work \cite{Jiang2024}, where the backward differentiation formula is utilized to discretize the time variable, resulting in the following scheme, referred to here as the BGN/BDF2 scheme:
	
\noindent(\textbf{Second-order BGN/BDF2 scheme}):~ Given $\bX^0,\bX^1 \in [V^h]^2$  which are appropriate approximations for $\Gamma(0)$ and $\Gamma(\tau)$, respectively, find $\mathbf{X}^{m+1}\in [V^h]^2$ and $\kappa^{m+1}\in V^h$ for $m\ge 1$ such that
\begin{equation}\label{SDF:BDF2}
\begin{split}
\l(\frac{\frac{3}{2}\mathbf{X}^{m+1}-2\mathbf{X}^{m}+\frac{1}{2}\mathbf{X}^{m-1}}{\tau},\varphi^h \widetilde{\mathbf{n}}^{m+1} \r)^h_{\widetilde{\Gamma}^{m+1}}&+\l(\p_s\kappa^{m+1},\p_s\varphi^h \r)_{\widetilde{\Gamma}^{m+1}}=0,\qquad \forall\ \varphi^h\in V^h,\\
\l(\kappa^{m+1},\widetilde{\mathbf{n}}^{m+1}\cdot \bm{\omega}^h\r)_{\widetilde{\Gamma}^{m+1}}^h&-\l(\p_s\mathbf{X}^{m+1},
\p_s\bm{\omega}^h\r)_{\widetilde{\Gamma}^{m+1}}=0,\quad \,\,\, \forall\ \bm{\omega}^h\in  [V^h]^2.
\end{split}
\end{equation}
Here  $\widetilde{\Gamma}^{m+1}$, described by $\widetilde{\bX}^{m+1}\in [V^h]^2$, provides a suitable approximation of
$\Gamma(t_{m+1})$ as predicted by the classical BGN scheme with a time step of $\tau$. Both the normal vector $\widetilde{\mathbf{n}}^{m+1}$ and the derivative $\p_s$ are defined over the polygon $\widetilde{\Gamma}^{m+1}$. The derived scheme also exhibits second-order accuracy. However, numerical results indicate that the BGN/PC scheme achieves a more robust convergence order and superior accuracy compared to the BGN/BDF2 scheme (see Figures \ref{Fig:EOC} and \ref{Fig:SDF3d_EOC}).
\end{remark}

\smallskip

We subsequently establish several fundamental properties of the BGN/PC scheme \eqref{SDF:BGN/PC}. We first demonstrate the well-posedness of BGN/PC scheme \eqref{SDF:BGN/PC} under certain mild assumptions regarding the predicted curve $\widetilde{\Gamma}^{m+\frac{1}{2}}$.

\begin{thm}[Well-posedness]\label{Well-posedness}
For $m\ge 1$, we assume that the discrete polygon $\widetilde{\Gamma}^{m+\frac{1}{2}}:=\widetilde{\bX}^{m+\frac{1}{2}}$ satisfies the following conditions within the BGN/PC algorithm:
\begin{enumerate}
		\item[(1)] At least two vectors from $\{\widetilde{\mathbf{h}}^{m+\frac{1}{2}}_j\}_{j=1}^{N}$ are not parallel, i.e.,
		\[
		\mathrm{dim}\l( \mathrm{Span}\l\{\widetilde{\mathbf{h}}^{m+\frac{1}{2}}_j \r\}_{j=1}^{N}\r)=2;
		\]
		\item[(2)] There are no degenerate elements on $\widetilde{\Gamma}^{m+\frac{1}{2}}$; that is,
		$\min_{1\le j\le N}|\widetilde{ \mathbf{h}}^{m+\frac{1}{2}}_j|>0$.
	\end{enumerate}
Then, the BGN/PC scheme \eqref{SDF:BGN/PC} is well-posed, i.e., there exists a unique solution $(\mathbf{X}^{m+1},\kappa^{m+1})\in [V^h]^2\times V^h$ to \eqref{SDF:BGN/PC}.
	
\end{thm}

\begin{proof}
	The proof follows as a specific instance of the anisotropic case by setting $\gamma(\bn)\equiv 1$ and $\bZ_k(\bn)=I_2$. For a detailed demonstration, we refer to Theorem \ref{aSDF:PC}.
\end{proof}

\smallbreak

In contrast  to the BGN/CNLF scheme \eqref{SDF:BGN2}, we can prove that our BGN/PC scheme \eqref{SDF:BGN/PC} retains the long-time equidistribution property characteristic of the BGN scheme \eqref{SDF:BGN1}. The proof is inspired by the idea in \cite[Proposition 3.3]{Zhao2021} concerning the first-order scheme for solving the isotropic surface diffusion of open curves.

\smallskip

\begin{thm}\label{BGN/PC:Equidistributed property}
	 Let $(\bX^{m},\kappa^m)$ be the solution of the BGN/PC scheme \eqref{SDF:BGN/PC}. We assume the following conditions:
	\begin{enumerate}
	    \item[(1)] As $m\rightarrow +\infty$, $\bX^{m}$ converges to the equilibrium state $\Gamma^e:=\bX^e\in [V^h]^2$, and $\kappa^m$ converges to $\kappa^e\in V^h$. Moreover, we assume the predictor $\widetilde{\bX}^{m+\frac{1}{2}}$ also converges to the same equilibrium $\Gamma^e$;
		
		\item[(2)] The  equilibrium polygon $\Gamma^e$  is  non-degenerate in the sense that
 \[\min_{1\le j\le N}|\mathbf{h}_j^e|>0,\quad  \text{and}\quad \mathrm{dim}\l( \mathrm{Span}\l\{\mathbf{h}^{e}_j \r\}_{j=1}^{N}\r)=2,\quad \mathbf{h}_j^e:=\bX^e(\rho_j)-\bX^e(\rho_{j-1}),\quad 1\le j\le N.\]
	\end{enumerate}
Then, the mesh ratio $\Psi^m:=\frac{\max_{1\le j\le N}|\mathbf{h}_j^m|}{\min_{1\le j\le N}|\mathbf{h}_j^m|}$ satisfies
\[
	\lim_{m\rightarrow +\infty}\Psi^m=\Psi^e:=\frac{\max_{1\le j\le N}|\mathbf{h}_j^e|}{\min_{1\le j\le N}|\mathbf{h}_j^e|}=1.
\]
Moreover, the equilibrium state $\Gamma^e$ is a regular $N$-polygon.
\end{thm}

\begin{proof}
  By taking $m\rightarrow +\infty$ in the BGN/PC scheme \eqref{SDF:BGN/PC}, alongside assumption (1), we derive
\begin{equation}\label{corrector:equi}
		\begin{split}
	&\l( \p_s\kappa^e,\p_s\varphi^h \r)_{\Gamma^{e}}=0,\quad \forall\ \varphi^h\in V^h,\\
	&\l(\kappa^e,\mathbf{n}^{e}\cdot \bm{\omega}^h\r)^h_{\Gamma^{e}}-\l(\p_s\mathbf{X}^{e},\p_s\bm{\omega}^h\r)_{\Gamma^{e}}=0, \quad \forall\ \bm{\omega}^h\in  [V^h]^2.
		\end{split}
	\end{equation}
Taking $\varphi^h=\kappa^e\in V^h$ in the first equation yields that
$\l( \p_s\kappa^e,\p_s\kappa^e\r)_{\Gamma^{e}}=0$. Noticing $\kappa^e\in \mathbb{C}(\mathbb{I})$ is a constant, this implies $\kappa^e\equiv\kappa^c\in\mathbb{R}$. Through direct computation, we can reformulate the second equation in \eqref{corrector:equi} into the following vector equation
\begin{equation}\label{Vectorized equation} \kappa^c(\mathbf{h}^e_{j+1}+\mathbf{h}^e_{j})^\perp+\l(\frac{\mathbf{h}^e_j}{|\mathbf{h}^e_j|}-\frac{\mathbf{h}^e_{j+1}}{|\mathbf{h}^e_{j+1}|} \r)=0, \quad j=1,\ldots,N.
\end{equation}
Multiplying both sides of \eqref{Vectorized equation} by the vector $(\mathbf{h}^e_{j+1}+\mathbf{h}^e_{j})$, we obtain
\begin{align*}
	0 =\l(\frac{\mathbf{h}^e_j}{|\mathbf{h}^e_j|}-\frac{\mathbf{h}^e_{j+1}}{|\mathbf{h}^e_{j+1}|} \r)\cdot (\mathbf{h}^e_{j+1}+\mathbf{h}^e_{j})=\l(|\mathbf{h}^e_{j+1}|-|\mathbf{h}^e_{j}|\r)\l(\frac{\mathbf{h}^e_{j+1}\cdot \mathbf{h}^e_{j}}{|\mathbf{h}^e_{j+1}||\mathbf{h}^e_{j}|}-1\r).
\end{align*}
Applying Cauchy-Schwarz inequality, we conclude that, for $j=1,\ldots,N$,
\[
|\mathbf{h}^e_{j+1}|=|\mathbf{h}^e_{j}|,\quad \text{or}\quad \mathbf{h}^e_{j+1}\parallel \mathbf{h}^e_{j}.
\]
Note that if the latter scenario occurs for some $j$, then $\kappa^c=0$ and all vectors $\mathbf{h}^e_{j}$ must be parallel, which contradicts assumption (2). Thus, we establish the equidistribution property for $\Gamma^e$, namely,
\[
|\mathbf{h}^e_{j}|\equiv |\mathbf{h}^e|,\quad j=1,\ldots,N.
\]

It remains to demonstrate all the interior angles are equal. By multiplying both sides of \eqref{Vectorized equation} by the vector $(\mathbf{h}^e_{j+1}+\mathbf{h}^e_{j})^{\perp}$  and denoting $\alpha_j\in (-\pi,\pi]$ as the exterior angle formed by $\mathbf{h}^e_{j}$ and $\mathbf{h}^e_{j+1}$, which is uniquely determined by
\[
\mathbf{h}^e_{j+1}\cdot \mathbf{h}^e_{j}=|\mathbf{h}^e_{j+1}||\mathbf{h}^e_{j}|\cos\alpha_j,\quad\text{and}\quad \mathbf{h}^e_{j}\cdot \l(\mathbf{h}^e_{j+1}\r)^\perp=-|\mathbf{h}^e_{j+1}||\mathbf{h}^e_{j}|\sin\alpha_j,
\]
we then obtain
\begin{align*}
	0 &=\kappa^c|\mathbf{h}^e_{j+1}+\mathbf{h}^e_{j}|^2+\l(\frac{\mathbf{h}^e_{j}}{|\mathbf{h}^e_{j}|}-\frac{\mathbf{h}^e_{j+1}}{|\mathbf{h}^e_{j+1}|} \r)\cdot \l(\mathbf{h}^e_{j+1}+\mathbf{h}^e_{j}\r)^\perp\\	&=\kappa^c|\mathbf{h}^e_{j+1}+\mathbf{h}^e_{j}|^2+\l(\frac{|\mathbf{h}^e_{j+1}|+|\mathbf{h}^e_{j}|}{|\mathbf{h}^e_{j+1}||\mathbf{h}^e_{j}|} \r)\mathbf{h}^e_{j}\cdot \l(\mathbf{h}^e_{j+1}\r)^\perp\\ &=2|\mathbf{h}^e|^2\kappa^c\l(1+\cos\alpha_j\r)-2|\mathbf{h}^e|\sin\alpha_j.
\end{align*}
Based on the sign of the constant $\kappa^c$, we can identify three cases:
\begin{itemize}
	\item If $\kappa^c=0$, then \eqref{Vectorized equation} implies that all vectors $\mathbf{h}^e_j$ are linearly dependent, contradicting assumption (2);
	\item If $\kappa^c<0$, then $-\pi<\alpha_j<0$ for all $j$, which contradicts the definition of a polygon;
	\item If $\kappa^c>0$, then $0<\alpha_j< \pi$, allowing us to represent the constant $\kappa^c$ as $	\kappa^c=\frac{\sin\alpha_j}{(1+\cos\alpha_j)|\mathbf{h}^e|}$.
	Moreover, we conclude that all exterior angles $\alpha_j$ are equal since the function $f(\alpha)=\frac{\sin\alpha}{1+\cos\alpha}$ is strictly monotonic in the interval $0<\alpha< \pi$.
\end{itemize}
In summary, the mesh points of the equilibrium state $\Gamma^e$ are perfectly equidistributed, and it is a regular $N$-polygon. Therefore, we have completed the proof.
\end{proof}

\smallskip

\section{Extensions to other geometric flows}

In this section, we extend the predictor-corrector strategy to other geometric flows, including CSF, AP-CSF, A-SDF and SDF in $\R^3$. The basic idea closely resembles the BGN/PC algorithm \eqref{SDF:BGN/PC},
wherein we first predict the solution at an intermediate time level using a first-order method. We then apply a formally second-order method to refine and achieve a more accurate corrected solution.

\subsection{For curve shortening flow}

Based on the variational formulation \eqref{CSF/AP-CSF:Coupled equation}, the first-order scheme, as outlined in \cite{Barrett2020,Jiang2023}, involves replacing the first equation in \eqref{SDF:BGN1} with the following equations
\begin{equation}\label{CSFbgn1}
	\l(\frac{\mathbf{X}^{m+1}-\mathbf{X}^m}{\tau},\varphi^h \mathbf{n}^m \r)^h_{\Gamma^m}=-\begin{cases}
		\l(  \kappa^{m+1},\varphi^h \r)_{\Gamma^m}^h, &\qquad\qquad \text{CSF},\\
		\l(\kappa^{m+1}-\l<\kappa^{m+1}\r>_{\Gamma^m}^h ,\varphi^h \r)_{\Gamma^m}^h, & \qquad\qquad \text{AP-CSF},\\
	\end{cases}
\end{equation}
where $\l<\cdot\r>_{\Gamma^{m}}=\l(\cdot,1 \r)_{\Gamma^{m}}^h/\l(1,1 \r)_{\Gamma^{m}}^h $.

\smallskip
To extend the BGN/PC scheme \eqref{SDF:BGN/PC} for solving CSF/AP-CSF, we only need to substitute the first equation in \eqref{SDF:BGN/PC} with
	\begin{align}\label{CSFandAP-CSF}
		\l(\frac{\mathbf{X}^{m+1}-\mathbf{X}^{m}}{\tau},\varphi^h \widetilde{\mathbf{n}}^{m+\frac{1}{2}} \r)^h_{\widetilde{\Gamma}^{m+\frac{1}{2}}}=-\begin{cases}
			\l( \frac{\kappa^{m+1}+\kappa^{m}}{2},\varphi^h \r)_{\widetilde{\Gamma}^{m+\frac{1}{2}}}^h,\qquad\qquad\qquad\qquad\qquad\quad \text{CSF},\\
			\l( \frac{\kappa^{m+1}+\kappa^{m}}{2}- \l<\frac{\kappa^{m+1}+\kappa^{m}}{2}\r>_{\widetilde{\Gamma}^{m+\frac{1}{2}}},\varphi^h \r)_{\widetilde{\Gamma}^{m+\frac{1}{2}}}^h,\qquad \text{AP-CSF},
		\end{cases}
	\end{align}
where $\widetilde{\Gamma}^{m+\frac{1}{2}}$ and the normal vector $\widetilde{\bn}^{m+\frac{1}{2}}$ are obtained in a manner analogous to that used for the isotropic planar SDF. Specifically,  $\widetilde{\Gamma}^{m+\frac{1}{2}}$ represents the solution of \eqref{CSFbgn1} with initial data $\mathbf{X}^m$ and a time step of $\tau/2$, and $\widetilde{\bn}^{m+\frac{1}{2}}$ is derived as the normal vector of $\widetilde{\Gamma}^{m+\frac{1}{2}}$.
Additionally,  $\l<\cdot\r>_{\widetilde{\Gamma}^{m+\frac{1}{2}}}=\l(\cdot,1 \r)_{\widetilde{\Gamma}^{m+\frac{1}{2}}}^h/\l(1,1 \r)_{\widetilde{\Gamma}^{m+\frac{1}{2}}}^h $ is computed over the predicted polygon curve $\widetilde{\Gamma}^{m+\frac{1}{2}}$. The comprehensive algorithm follows a framework similar to the BGN/PC Algorithm, and for the sake of brevity, we omit the full details here.

\subsection{For anisotropic surface diffusion flow of planar curves}

In their seminal work \cite{Bao2022a}, Bao, Jiang and Li (referred to as BJL) proposed an energy-stable PFEM for solving anisotropic SDF of curves characterized by arbitrary surface energy densities $\gamma(\bn)$, which encompasses both weakly and strongly anisotropic surface energy cases. Specifically, they reformulated the coupled equations \eqref{aSDF:Coupled equation2} using the Cahn-Hoffman $\bm{\xi}$-vector formulation alongside a vital identity \cite[Lemma 2.1]{Bao2022a}
\[
\bZ_k(\bn)\p_s\bX=\bm{\xi}^\perp,
\]
where  $\bZ_k(\bn)$ is defined as
\begin{equation}\label{Def of Z}
\bZ_k(\bn)=\gamma(\bn)I_2-\bn\bm{\xi}^{\top}-\bm{\xi}\bn^{\top}+k(\bn)\bn\bn^{\top},\quad \forall\, \bn\in \S^1,
\end{equation}
and $k(\bn):\S^1\rightarrow \R^+$ serves as a stabilizing function to ensure that  $\bZ_k(\bn)$ is a positive definite matrix. Here $I_2$ represents the identity matrix, and $\cdot^\top$ denotes the transpose operator. Utilizing this variational formulation, BJL proposed the following first-order BJL scheme:

\noindent(\textbf{First-order BJL scheme}):~For $m\ge 0$, find $\mathbf{X}^{m+1}\in [V^h]^2$ and $\mu^{m+1}\in V^h$ such that
\begin{equation}\label{aSDF:fisrt-order}
		\begin{cases}
			\l(\frac{\mathbf{X}^{m+1}-\mathbf{X}^m}{\tau},\varphi^h \mathbf{n}^m \r)^h_{\Gamma^m}+\l( \p_s \mu^{m+1},\p_s \varphi^h \r)_{\Gamma^m}=0,\quad \forall\ \varphi^h\in V^h,\\
			\vspace{-3mm}\\
				\l(\mu^{m+1},\mathbf{n}^m\cdot \bm{\omega}^h\r)_{\Gamma^m}^h-\l(\bZ_k(\bn^m) \p_s \mathbf{X}^{m+1},\p_s\bm{\omega}^h\r)_{\Gamma^m}=0,\quad \forall\ \bm{\omega}^h\in  [V^h]^2.
		\end{cases}
	\end{equation}
 The well-posedness and unconditional energy stability have be established under mild conditions. It should be noted that the existence of the stabilizing function $k(\bn)$ relies on specific assumptions pertaining to anisotropic surface energy densities $\gamma(\bn)$, such as
 \begin{equation}\label{gamma assumption}
 	\gamma(-\bn)=\gamma(\bn),\quad \forall\, \bn\in \S^1,\quad \gamma(\mathbf{p})\in C^2(\R^2\setminus\{\mathbf{0}\}).
 \end{equation}
 For explicit computation of the stabilizing function $k(\bn)$, we refer readers to \cite[Appendix]{Bao2022a}. In practice, the scheme \eqref{aSDF:fisrt-order} demonstrates quadratic convergence in space \cite{Bao2022a} and linear convergence in time.

The predictor-corrector strategy outlined above can be readily extended to  the anisotropic case using  the first-order BJL  scheme \eqref{aSDF:fisrt-order}. Analogously, the terms $\p_t\bX$, $\mu$ and $\p_s \bX$ are approximated using the Taylor expansion formulas \eqref{PC method approximation}. Therefore, we propose the following second-order BJL/PC scheme:

\smallskip

(\textbf{Second-order BJL/PC scheme}):~ For $\bX^0 \in [V^h]^2$, $\mu^0\in  V^h$  which are appropriate approximations at $t_0=0$. For $m\ge 0$, let $\widetilde{\Gamma}^{m+\frac{1}{2}}:=\widetilde{\bX}^{m+\frac{1}{2}}\in [V^h]^2$ be an appropriate approximation of $\bX(t_{m+\frac{1}{2}})$; we seek $\mathbf{X}^{m+1}\in [V^h]^2$ and $\mu^{m+1}\in V^h$  such that
\begin{equation}\label{aSDF:PC}
		\begin{cases}
	\l(\frac{\mathbf{X}^{m+1}-\mathbf{X}^{m}}{\tau},\varphi^h \widetilde{\mathbf{n}}^{m+\frac{1}{2}} \r)^h_{\widetilde{\Gamma}^{m+\frac{1}{2}}}+\l( \p_s\l(\frac{\mu^{m+1}+\mu^{m}}{2}\r),\p_s\varphi^h \r)_{\widetilde{\Gamma}^{m+\frac{1}{2}}}=0,\quad \forall\ \varphi^h\in V^h,\\
	\vspace{-3mm}\\		\l(\frac{\mu^{m+1}+\mu^{m}}{2},\widetilde{\mathbf{n}}^{m+\frac{1}{2}}\cdot \bm{\omega}^h\r)_{\widetilde{\Gamma}^{m+\frac{1}{2}}}^h-\l(\bZ_k( \widetilde{\mathbf{n}}^{m+\frac{1}{2}}) \p_s\l(\frac{ \mathbf{X}^{m+1}+ \mathbf{X}^{m}}{2}\r),\p_s\bm{\omega}^h\r)_{\widetilde{\Gamma}^{m+\frac{1}{2}}}=0, \quad \forall\ \bm{\omega}^h\in  [V^h]^2.
		\end{cases}
	\end{equation}
The computation of the stabilizing function $k(\bn)$ and $\bZ_{k}(\bn)$ is exactly the same as the first-order BJL scheme \eqref{aSDF:fisrt-order} (see \cite[Appendix]{Bao2022a}). Similar to the BGN/PC algorithm, we  propose the following second-order algorithm as follows.

\begin{algorithm}\textbf{BJL/PC algorithm for anisotropic case}
	\label{second-order algorithm for ani}
\begin{algorithmic}[1]
 \REQUIRE An initial curve $\Gamma(0)$ approximated by a polygon $\Gamma^0$ with $N$ vertices, described by $\bX^0\in [V^h]^2$, time step $\tau$ and terminate time $T$ satisfying $T/\tau\in\mathbb{N}$.
 \ENSURE Computational solution $\Gamma^{T/\tau}:=\bX^{T/\tau}$.
 \STATE Calculate $\kappa^0$ by  $\Gamma^0$ and a formula based on the formula $\kappa\bn=-\p_{ss}\bX$ and least square method \cite{Barrett2007a,Jiang2023}, and we further set the initial chemical potential by the $\theta$-formulation \cite{Bao17}
 \begin{equation}
 	\mu^0=(\gamma(\theta^0)+\gamma^{\prime\prime}(\theta^0))\, \kappa^0,
 \end{equation}
 where $\theta^0$ is the angle between $y$-axis and the normal vector $\bn^0$ of the initial polygon $\Gamma^0$. Set $m=0$.
 \WHILE{$m < T/\tau$,}
 \STATE {\small
 \begin{itemize}
\item \textbf{Predictor}: Compute  $\widetilde{\Gamma}^{m+\frac{1}{2}}:= \widetilde{\bX}^{m+\frac{1}{2}}$ by using the   BJL scheme \eqref{aSDF:fisrt-order} with $\Gamma^m$ and $\tau/2$.
 \item \textbf{Corrector}: Compute $\Gamma^{m+1}:=\bX^{m+1},\mu^{m+1}$ by using the   BJL/PC scheme  \eqref{aSDF:PC} with
 $\bX^{m},\mu^{m},\widetilde{\Gamma}^{m+\frac{1}{2}}$ and $\tau$.
\end{itemize}
}
 \STATE $m = m + 1$;
 \ENDWHILE
 \end{algorithmic}
\end{algorithm}

\begin{thm}[Well-posedness]\label{ani-Well-posedness}
For $m\ge 1$,  assume that the predicted polygon curve $\widetilde{\Gamma}^{m+\frac{1}{2}}:=\widetilde{\bX}^{m+\frac{1}{2}}$ satisfies conditions (1) and (2) as outlined in Theorem \ref{Well-posedness}. Additionally, we assume that the surface energy density satisfies \eqref{gamma assumption}. Under these assumptions, the BJL/PC  scheme \eqref{aSDF:PC} is well-posed, i.e., there exists a unique solution $(\mathbf{X}^{m+1},\mu^{m+1})\in [V^h]^2\times V^h$.
\end{thm}

\begin{proof}
	Notably, given the predicted polygon $\widetilde{\Gamma}^{m+\frac{1}{2}}$, the resulting scheme \eqref{aSDF:PC} generates a system of linear algebraic equations with respect to the unknowns $(\bX^{m+1},\mu^{m+1})$. To establish well-posedness, it suffices to demonstrate that the following homogeneous system for $(\bX,\mu)\in [V^h]^2\times V^h$ has only the trivial zero solution:
	\begin{equation}\label{Linear equation}
		\begin{cases}
			\l(\frac{\mathbf{X}}{\tau},\varphi^h \widetilde{\mathbf{n}}^{m+\frac{1}{2}} \r)^h_{\widetilde{\Gamma}^{m+\frac{1}{2}}}+\l( \p_s\mu ,\p_s\varphi^h \r)_{\widetilde{\Gamma}^{m+\frac{1}{2}}}=0,\quad \forall\ \varphi^h\in V^h,\\
			\vspace{-3mm}\\
   \l(\mu,\widetilde{\mathbf{n}}^{m+\frac{1}{2}}\cdot \bm{\omega}^h\r)_{\widetilde{\Gamma}^{m+\frac{1}{2}}}^h-\l(\bZ_k( \widetilde{\mathbf{n}}^{m+\frac{1}{2}})\p_s \mathbf{X},\p_s\bm{\omega}^h\r)_{\widetilde{\Gamma}^{m+\frac{1}{2}}}=0,\quad \forall\ \bm{\omega}^h\in  [V^h]^2.
		\end{cases}
	\end{equation}
By letting $\varphi^h=\mu $ and $\bm{\omega}^h=\bX$, we can derive that
\begin{equation}
	 \l( \p_s\mu ,\p_s\mu \r)_{\widetilde{\Gamma}^{m+\frac{1}{2}}}+\l(\bZ_k( \widetilde{\mathbf{n}}^{m+\frac{1}{2}})\p_s \mathbf{X},\p_s\bX\r)_{\widetilde{\Gamma}^{m+\frac{1}{2}}}=0.
\end{equation}
 Given the assumption \eqref{gamma assumption}, the choice of the stabilizing function $k(\bn)$ ensures that $\bZ_k(\widetilde{\mathbf{n}}^{m+\frac{1}{2}})$ is a symmetric positive definite matrix. Consequently,  we conclude that $\bX$ and $\mu$ must be constants:
\begin{equation*}
		\bX\equiv \bX^c\in \R^2,\quad \mu \equiv \mu^c\in \R.
	\end{equation*}
Substituting this into \eqref{Linear equation} leads to
\begin{equation}
	\l(\mathbf{X}^c,\varphi^h \widetilde{\mathbf{n}}^{m+\frac{1}{2}} \r)^h_{\widetilde{\Gamma}^{m+\frac{1}{2}}}=0,\quad   \l(\mu^c,\widetilde{\mathbf{n}}^{m+\frac{1}{2}}\cdot \bm{\omega}^h\r)_{\widetilde{\Gamma}^{m+\frac{1}{2}}}^h=0,\quad \forall\ (\varphi^h,\bm{\omega}^h)\in V^h\times[V^h]^2.
\end{equation}
The standard argument presented in \cite[Theorem 2.1]{Barrett2007a} leads to the conclusion that $\bX^c=0$ and $\mu^c=0$ under conditions (1) and (2) for $\widetilde{\Gamma}^{m+\frac{1}{2}}$. This completes the proof of well-posedness.
\end{proof}

\subsection{For surface diffusion flow of surfaces in $\R^3$}

In alignment with the previous discussion and utilizing the first-order BGN scheme for solving SDF in $\R^3$ as outlined in \cite{Barrett2008c}, we extend the predictor-corrector strategy to derive a second-order scheme in this subsection.

We approximate the surface $\Gamma(t_m)$ by a polyhedron $\Gamma^m:=\bigcup^J_{j=1} \overline{\sigma^m_j}$, where $\{\sigma^m_j\}_{j=1}^J$ constitutes a family of mutually disjoint open triangles. Set $h=\max\limits_{1\le j\le J}\mathrm{diam}(\sigma_j^m)$ and consider the finite element space $W_m^h$, which consists of piecewise linear functions defined on $\Gamma^m$:
\[
W_m^h:=\l\{u\in C(\Gamma^m, \mathbb{R}):\ u|_{\sigma^m_j}\ \text{is linear},\ \forall j=1,\ldots,J \r\}\subseteq H^1(\Gamma^m, \mathbb{R}).
\]
We apply the following mass-lumped inner product to approximate the $L^2$ inner product $( \cdot, \cdot)_{\Gamma^m}$:
\begin{align*}
\l(u,v\r)_{\Gamma^m}^h:=\frac{1}{3}\sum_{j=1}^J|\sigma^m_j|
\sum_{k=1}^3(u\cdot v)\l((\mathbf{q}_{j_k}^m)^{-} \r),
\end{align*}
where $\{\mathbf{q}_{j_1},\mathbf{q}_{j_2},\mathbf{q}_{j_3}\}$ are the vertices of the triangle $\sigma^m_j$ and $|\sigma^m_j|$ is its area. Additionally, $u\l((\mathbf{q}_{j_k}^m)^{-} \r)=\lim\limits_{\sigma^m_j\ni \mathbf{x}\rightarrow \mathbf{q}_{j_k}^m}u(\mathbf{x}).$
The first-order BGN scheme presented in \cite{Barrett2008c} for solving SDF of surfaces in $\R^3$ can be written as follows:

\smallskip

\noindent(\textbf{First-order BGN scheme}): Given $\Gamma^0$ and its parametrized identity function $\bX^0\in W_0^h$ on $\Gamma^0$, for $m\ge 0$, we seek $\mathbf{X}^{m+1}\in [W_m^h]^3$ and $H^{m+1}\in W_m^h$ such that
\begin{equation}\label{SDF3d:BGN1}
		\begin{cases}
			\l(\frac{\mathbf{X}^{m+1}-\mathbf{X}^m}{\tau},\varphi^h \mathbf{n}^m \r)^h_{\Gamma^m}+\l( \nabla_{\Gamma^m} H^{m+1},\nabla_{\Gamma^m}\varphi^h \r)_{\Gamma^m}=0,\quad \forall\ \varphi^h\in W_m^h,\\	\l(H^{m+1},\mathbf{n}^m\cdot \bm{\omega}^h\r)_{\Gamma^m}^h-\l(\nabla_{\Gamma^m}\mathbf{X}^{m+1},\nabla_{\Gamma^m} \bm{\omega}^h\r)_{\Gamma^m}=0,\quad \forall\ \bm{\omega}^h\in  [W_m^h]^3,
		\end{cases}
	\end{equation}
where $\mathbf{X}^m(\cdot)$ is the identity function on $[W_m^h]^3$, the unit outward normal vector $\mathbf{n}^m$ and the surface gradient $\nabla_{\Gamma^m}$ are both piecewisely defined over the polyhedra $\Gamma^m$.

The predictor-corrector strategy is now readily applied in this context, leading to the following second-order BGN/PC scheme for solving SDF of surfaces in $\R^3$:

\smallskip

\noindent(\textbf{Second-order BGN/PC scheme}):~ For $\bX^0 \in [W_0^h]^3$ and $H^0\in  W_0^h$, which are appropriate approximations at the time level $t_0=0$. For $m\ge 0$, let $\widetilde{\Gamma}^{m+\frac{1}{2}}:=\widetilde{\bX}^{m+\frac{1}{2}}\in [W_m^h]^3$ be an appropriate approximation of $\bX(t_{m+\frac{1}{2}})$. We seek $\mathbf{X}^{m+1}\in [W_m^h]^3$ and $H^{m+1}\in W_m^h$ such that
\begin{equation}\label{SDF3d:BGN/PC}
		\begin{cases}
			\l(\frac{\mathbf{X}^{m+1}-\mathbf{X}^m}{\tau},\varphi^h \widetilde{\mathbf{n}}^{m+\frac{1}{2}} \r)^h_{\widetilde{\Gamma}^{m+\frac{1}{2}}}+\l( \nabla_{\widetilde{\Gamma}^{m+\frac{1}{2}}}\l(\frac{H^{m+1}+H^{m}}{2}\r),\nabla_{\widetilde{\Gamma}^{m+\frac{1}{2}}}\varphi^h \r)_{\widetilde{\Gamma}^{m+\frac{1}{2}}}=0,\ \forall\ \varphi^h\in W_m^h,\\	\l(\frac{H^{m+1}+H^{m}}{2},\widetilde{\mathbf{n}}^{m+\frac{1}{2}}\cdot \bm{\omega}^h\r)_{\widetilde{\Gamma}^{m+\frac{1}{2}}}^h-\l(\nabla_{\widetilde{\Gamma}^{m+\frac{1}{2}}}\l(\frac{\mathbf{X}^{m+1}+\mathbf{X}^{m}}{2}\r),\nabla_{\widetilde{\Gamma}^{m+\frac{1}{2}}} \bm{\omega}^h\r)_{\widetilde{\Gamma}^{m+\frac{1}{2}}}=0,\ \forall\ \bm{\omega}^h\in  [W_m^h]^3.
		\end{cases}
	\end{equation}
In this setup, the unit normal vector $\widetilde{\mathbf{n}}^{m+\frac{1}{2}}$ and the surface gradient $\nabla_{\widetilde{\Gamma}^{m+\frac{1}{2}}}$ are both piecewisely defined over the predicted polyhedron $\widetilde{\Gamma}^{m+\frac{1}{2}}$, and they are obtained using the first-order scheme \eqref{SDF3d:BGN1} with a time step of $\tau/2$. The second-order BGN/PC Algorithm can be given similarly to the planar case; for conciseness, further details are omitted here.
\section{Numerical results}
In this section, we present extensive numerical experiments to demonstrate the efficiency and accuracy of our proposed BGN/PC algorithms for solving planar isotropic  SDF and other geometric flows.

As discussed in our previous work (e.g., see Section 3 in \cite{Jiang2023}), shape metrics are more suitable for quantifying numerical errors of BGN-based schemes, which introduce the intrinsic tangential velocities to maintain an even distribution of mesh points during evolution. To evaluate the convergence order, we employ the manifold distance to  measure the differences between two curves/surfaces. Specifically, the manifold distance between two curves/surfaces $\Gamma_1$ and $\Gamma_2$ is defined as~\cite{Zhao2021,Jiang2023}
\begin{align*}
	\mathrm{M}\l(\Gamma_1,\Gamma_2\r)
	&: = |(\Omega_1\setminus\Omega_2)\cup (\Omega_2\setminus\Omega_1) | =|\Omega_1 |+|\Omega_2 |-2 |\Omega_1\cap \Omega_2 |,
\end{align*}
where $\Omega_1$ and $\Omega_2$ denote the regions enclosed by $\Gamma_1$ and $\Gamma_2$, respectively, and $|\Omega|$ represents the area/volume of $\Omega$ in 2D/3D.
It can be easily proven that the manifold distance satisfies the properties of symmetry, positivity and the triangle inequality (see \cite[Proposition 5.1]{Zhao2021}), making it a reliable shape metric for measuring differences between curves in 2D or surfaces in 3D  \cite{Zhao2021,Bao2021,Bao2022b,Bao2023a,Jiang2023,Jiang2024}.

\subsection{For isotropic surface diffusion flow of planar curves}

\begin{ex}[Convergence order test]\label{ex:iso_EOC} In this experiment, taking the planar isotropic SDF as an example, we compare the convergence rates of four numerical schemes: the BGN scheme \eqref{SDF:BGN1}, the BGN/CNLF scheme \eqref{SDF:BGN2}, the BGN/BDF2 scheme \eqref{SDF:BDF2}, and  BGN/PC scheme \eqref{SDF:BGN/PC}, using an initially elliptic curve, parametrized as:
\[\bX(\rho)=(2\cos(2\pi \rho), \sin(2\pi \rho)),\quad \rho\in \mathbb{I}=[0,1].\]	
\end{ex}

Given the absence of an exact solution for the SDF corresponding to an initially elliptic curve, we compute a reference solution $\bX_{\mathrm{ref}}$ using the BGN/PC scheme with a very fine mesh and a tiny time step. Specifically, we choose $N=10^4$ and $\tau=10^{-1}\times 2^{-13}$. To assess the temporal convergence rates of the schemes, we maintain a large number of mesh points (e.g., $N=10^4$) in numerical simulations to ensure that the spatial error is negligible compared to the temporal error.
Subsequently, we calculate the numerical error and its corresponding convergence rate using the manifold distance, defined as follows:
\begin{equation}
 \label{eqn:errordef1}
\cE_{\tau}(T)=\cE_{M}(T)= \mathrm{M}(\bX^m_{\tau}, \bX_{\mathrm{ref}}),
\quad \text{Order}=\log\Big(\frac{\cE_{\tau_1}(T)}{\cE_{\tau_2}(T)} \Big)\Big/ \log \Big(\frac{\tau_1}{\tau_2} \Big),
\end{equation}
where $m=T/\tau$ is the discrete time level, and $\bX^m_{\tau}$ represents the numerical solution obtained from the various numerical schemes with time step $\tau$.

\begin{figure}[h!]
		\centering
		\includegraphics[width=15cm]{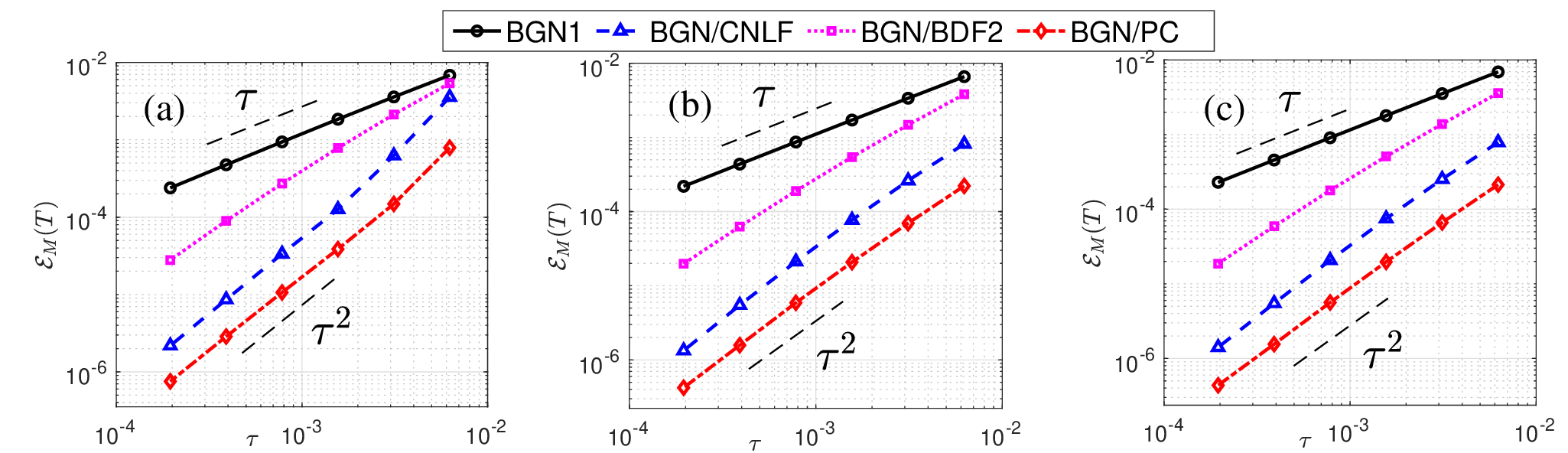}
		\setlength{\abovecaptionskip}{0pt}
		\vspace{10pt}
		\caption{Log-log plot of the numerical errors of the classical BGN scheme \eqref{SDF:BGN1}, the BGN/CNLF scheme \eqref{SDF:BGN2}, the BGN/BDF2 scheme \eqref{SDF:BDF2} and the BGN/PC scheme \eqref{SDF:BGN/PC}
		for solving the isotropic planar SDF associated with an initially elliptic curve at three time levels: (a) $T=0.05$, (b) $T=0.5$, (c) $T=5$.}
		\label{Fig:EOC}
	\end{figure}

As shown in Figure \ref{Fig:EOC}, we clearly observe that at three different time levels ($T=0.05, 0.5, 5$),
the numerical errors of the BGN/CNLF scheme \eqref{SDF:BGN2}, the BGN/BDF2 scheme \eqref{SDF:BDF2}, and the BGN/PC scheme \eqref{SDF:BGN/PC} all attain second-order convergence,
while the classical BGN scheme \eqref{SDF:BGN1} exhibits only first-order convergence.
Furthermore, it is noteworthy that the newly proposed BGN/PC scheme significantly outperforms the other three schemes in terms of  accuracy while utilizing the same computational parameters.

\smallskip

\begin{ex}[Comparison of computational cost] We compare the computational cost of the classical BGN scheme \eqref{SDF:BGN1}, the BGN/CNLF scheme \eqref{SDF:BGN2}, the BGN/BDF2 scheme \eqref{SDF:BDF2}, and the BGN/PC scheme \eqref{SDF:BGN/PC}. The experiments were performed using  MATLAB 2021b on a MacBook Pro with 1.4GHz quad-core Intel Core i5 processor and 8GB of RAM. The initial curve is still chosen as an elliptic curve, as defined previously.
\end{ex}

\begin{table}[htbp!]
	\centering
\caption{Comparisons of the CPU time (in seconds) and the numerical errors measured by the manifold distance $\cE_{M}(T)$ for the classical BGN scheme \eqref{SDF:BGN1}, the  BGN/CNLF scheme \eqref{SDF:BGN2}, the BGN/BDF2 scheme \eqref{SDF:BDF2} and the  BGN/PC scheme \eqref{SDF:BGN/PC}, where the initial curve is an ellipse, with computational parameters specified as $\tau=1/N$ and $T=0.05$.}
	\begin{tabular}{|c|c|c|c|c|c|c|c|c|}
\hline
\multicolumn{1}{|c|}{Mesh points } & \multicolumn{2}{|c|}{BGN  }  & \multicolumn{2}{|c|}{BGN/CNLF  } & \multicolumn{2}{|c|}{BGN/BDF2  }   & \multicolumn{2}{|c|}{BGN/PC}  \\ \hline
	   $N$  &  $\cE_{M}(T)$  &  \text{Time (s)} & $\cE_{M}(T)$  &   \text{Time (s)} & $\cE_{M}(T)$  &   \text{Time (s)}  & $\cE_{M}(T)$  &   \text{Time (s)}      \\ \hline
	   160  & 8.58E-3   & 0.14    & 3.65E-3    & 0.38  & 6.68E-3 & 0.25  & 1.79E-3   &  0.19     \\ \hline
	 320  &   4.02E-3  &  0.33    &  5.92E-4  &   0.45 & 2.45E-3 & 0.53 &  4.25E-4  & 0.77    \\ \hline
       640   &  1.96E-3  & 1.24     &   1.20E-4  &  1.86 & 8.67E-4 & 2.10  &  1.02E-4   &  3.10   \\ \hline
	1280    &  9.67E-4  & 6.86     &  2.94E-5 & 9.33 & 2.93E-4 & 12.8  & 2.47E-5  & 15.4  \\ \hline
      2560  &  4.81E-4  & 65.4    &  7.21E-6   &  87.3 & 9.46E-5 & 136 &  5.77E-6 &  154   \\ \hline
      5120   & 2.40E-4   & 1.19E+3   &  1.80E-6   & 1.37E+3  & 2.89E-5 & 2.24E+3  & 1.20E-6  & 2.57E+3 \\ \hline
\end{tabular}
\label{tab:CPU time:semi-implicit}
\end{table}

Table~\ref{tab:CPU time:semi-implicit} displays a comparison of CPU times (in seconds) and numerical errors at $T=0.05$, measured by the manifold distance $\cE_{M}(T)$, for the four numerical schemes mentioned above. The results indicate the following conclusions: (1) under the same computational parameters (i.e., $N$ and $\tau$), the numerical error of the BGN/PC scheme consistently outperforms that of the other three schemes;
(2) the computational cost of the BGN/CNLF scheme is comparable to that of the classical BGN scheme and is approximately half of the cost associated with the BGN/PC scheme and the BGN/BDF2 scheme. This is due to the additional calculations required for determining the predictor curve $\widetilde{\Gamma}^{m+\frac{1}{2}}$ at each time step for the BGN/PC scheme, which is similar to the the BGN/BDF2 scheme;
(3) when considering both accuracy and efficiency, the three second-order schemes demonstrate significantly better performance compared to the first-order BGN scheme \eqref{SDF:BGN1}.
Furthermore, taking all factors into consideration, the BGN/PC scheme performs the best among them.

\smallskip

\begin{ex}[Evolution of geometric quantities]\label{Evolution of geometric quantities}
We compare the evolution of several important geometric quantities using the mentioned four different schemes. The initial curve is still chosen as an elliptic curve. As is well-known, the SDF decreases the perimeter while preserving the enclosed area during evolution. This property inspires us to investigate the evolution of several geometric quantities: the normalized perimeter $L(t)/L(0)$, the relative area loss $\Delta A(t)$, and the mesh ratio function $\Psi(t)$, defined as follows:
\begin{align*}
    \l.\frac{L(t)}{L(0)}\r|_{t=t_m} = \frac{L^m}{L^0},\quad \l.\Delta A(t)\r|_{t=t_m} = \frac{A^m-A^0}{A^0},&\quad \l.\Psi(t)\r|_{t=t_m}=\frac{\max_{1\le j\le N}|\mathbf{h}_j^m|}{\min_{1\le j\le N}|\mathbf{h}_j^m|},\quad m\ge 0,
\end{align*}
where $L^m$ and $A^m$ represent the perimeter and the enclosed area of the polygon determined by the numerical solution $\bX^m$, respectively. The mesh ratio is closely related the mesh quality.
\end{ex}

\begin{figure}[htbp!]
		\centering
		\includegraphics[width=15cm]{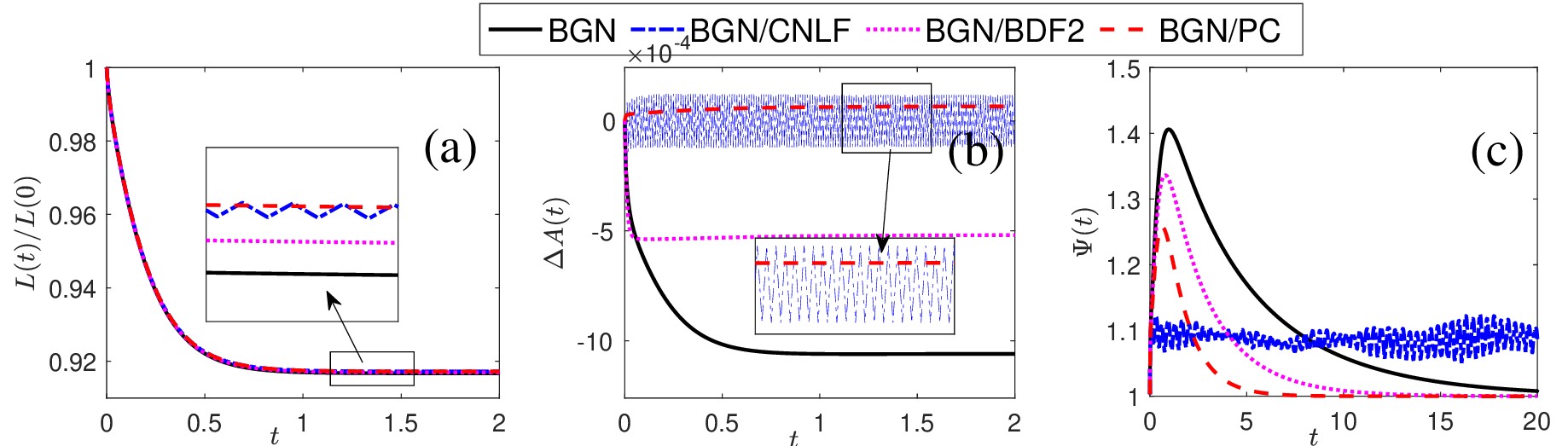}
			\setlength{\abovecaptionskip}{10pt}
		\caption{Evolution of several geometric quantities using the classical BGN scheme \eqref{SDF:BGN1}, the BGN/CNLF scheme \eqref{SDF:BGN2}, the BGN/BDF2 scheme \eqref{SDF:BDF2} and the BGN/PC  scheme \eqref{SDF:BGN/PC}:
		(a) normalized perimeter; (b) relative area loss; (c) mesh ratio, where the initial curve is an ellipse, with discretization parameters set as $N=80$ and $\tau=1/160$.}
		\label{Fig:GEO}
	\end{figure}

As illustrated in Figure  \ref{Fig:GEO}, we observe the following:  (1) all four numerical schemes generally decrease the perimeter of the polygon during the evolution (cf. Figure \ref{Fig:GEO}(a)). However, upon closer inspection, the BGN/CNLF scheme exhibits  oscillation behavior, likely due to its use of the two-step leap-frog method for temporal discretization;
(2) under identical computational parameters, the area loss of the BGN/PC scheme attains the least among the four schemes (cf. Figure \ref{Fig:GEO}(b)), further demonstrating its superior computational accuracy; (3) during long-time evolution, the mesh ratio function $\Psi(t)$ converges to $1$ for the BGN scheme, the BGN/BDF2 scheme, and the BGN/PC scheme (cf. Figure \ref{Fig:GEO}(c)), indicating asymptotic mesh equidistribution. In contrast, the mesh ratio for the BGN/CNLF scheme does not converge and again exhibits oscillatory behavior. In fact, this oscillation may lead to mesh distortion and potential breakdown of the BGN/CNLF scheme in certain scenarios, necessitating mesh regularization \cite{Jiang2023}.
Furthermore, it is noteworthy that among the three schemes which enjoy the mesh equidistribution property,
the BGN/PC scheme achieves this state most rapidly.

\begin{ex}[Simulation of morphological evolution]
	We employ the BGN/PC scheme \eqref{SDF:BGN/PC} to simulate the morphological evolution of various initially distinct curves towards the equilibrium state. Meanwhile, we illustrate the corresponding evolution of the geometric quantities discussed in the last example.
\end{ex}

\begin{figure}[h!]
		\centering
		\includegraphics[width=15cm,height=11cm]{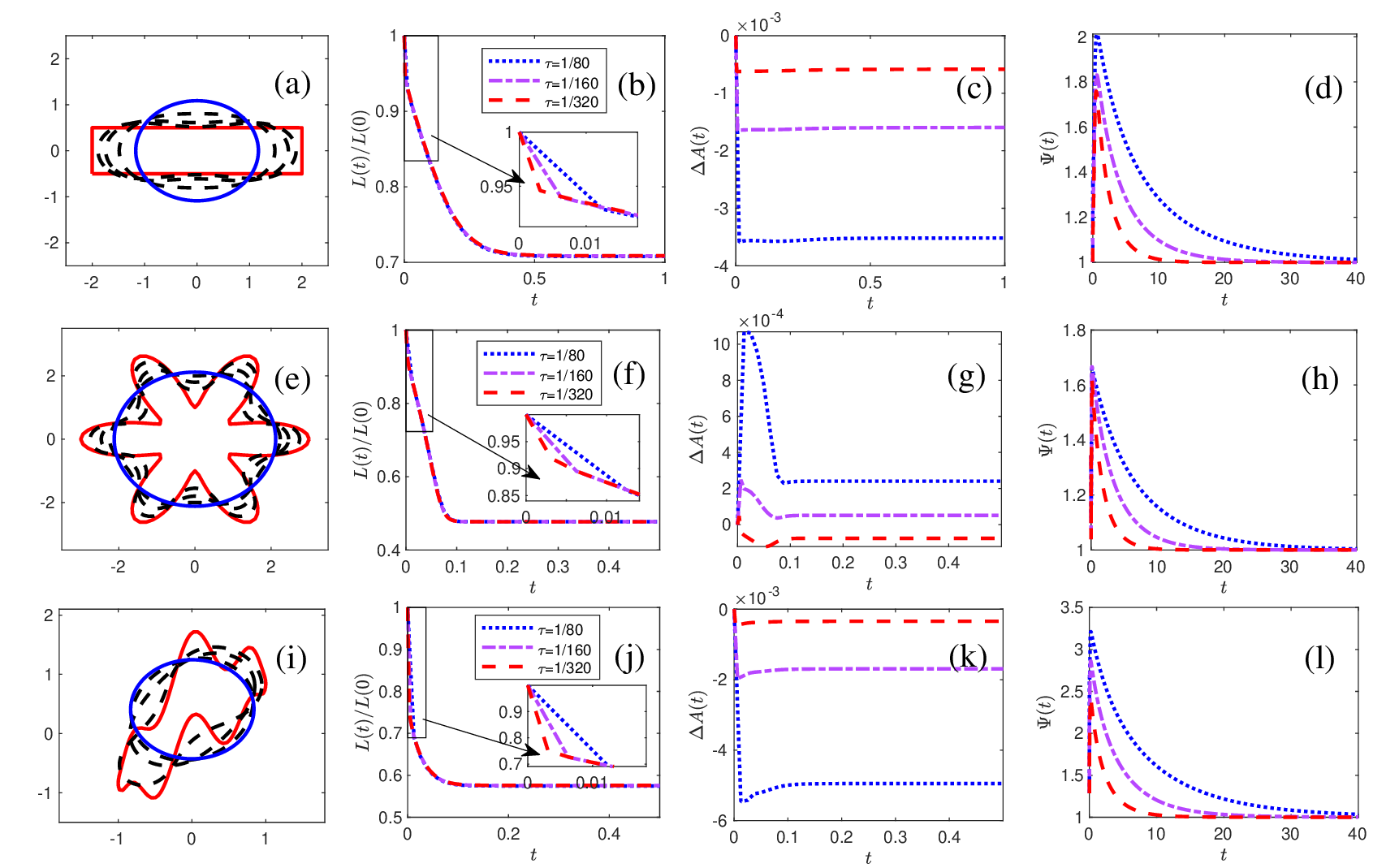}					\setlength{\abovecaptionskip}{5pt}
		\caption{Evolution of the curves and the corresponding geometric quantities by using the BGN/PC scheme \eqref{SDF:BGN/PC} for three different initial curves:
		 a rectangular curve (first row); a `flower' curve (second row); a non-convex curve (third row), where the discretization parameters are set as $N=320$, $\tau=1/80,1/160,1/320$.}
		\label{Fig:EVO_iso}
	\end{figure}

As illustrated in Figure \ref{Fig:EVO_iso}(a), (e) and (i), it is evident that regardless of the selected initial curve type,
all curves eventually evolve into the same equilibrium shape--a circle. Furthermore, by observing the evolution of geometric quantities, an examination of the evolution of geometric quantities leads us to the following conclusions:
(1) the BGN/PC scheme \eqref{SDF:BGN/PC} consistently decreases the perimeter during the evolution, aligning with the perimeter-decreasing property of the SDF (see  (b), (f) and (j)); (2) the BGN/PC scheme enjoys a long-time equidistribution property, whereby a smaller time step
$\tau$ accelerates the convergence of mesh points to the equidistribution state, i.e., the mesh ratio function
$\Psi(t)$ reaches $1$ more rapidly (cf. (d), (h) and (l)).

\medskip

\subsection{For curve shortening flow}

In this subsection, we apply the proposed BGN/PC scheme \eqref{CSFandAP-CSF} to the CSF  and AP-CSF, and present numerical results for both flows.

\begin{ex}[Extensions to CSF and AP-CSF]
We conduct convergence order tests of the BGN/PC scheme for solving CSF and AP-CSF, and illustrate the corresponding evolution of geometric quantities.
\end{ex}

\begin{figure}[h!]
		\centering
		\includegraphics[width=14cm,height=4cm]{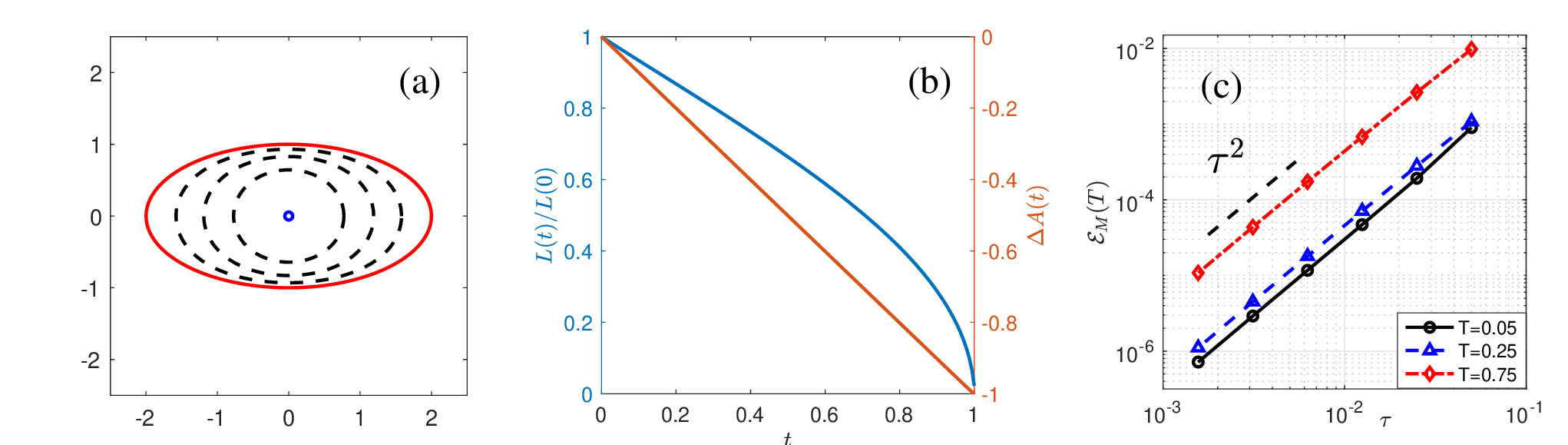}
		\caption{Evolution of an elliptic curve driven by the CSF using the BGN/PC scheme: (a) snapshots of the evolution; (b) evolution of the normalized perimeter and relative area loss; (c) temporal accuracy of the scheme at different time instances $T=0.05,0.25,0.75$, with the spatial mesh size fixed as $N=10^4$.}
		\label{Fig:CSF}
	\end{figure}

	\begin{figure}[h!]
		\centering
		\includegraphics[width=15cm,height=4cm]{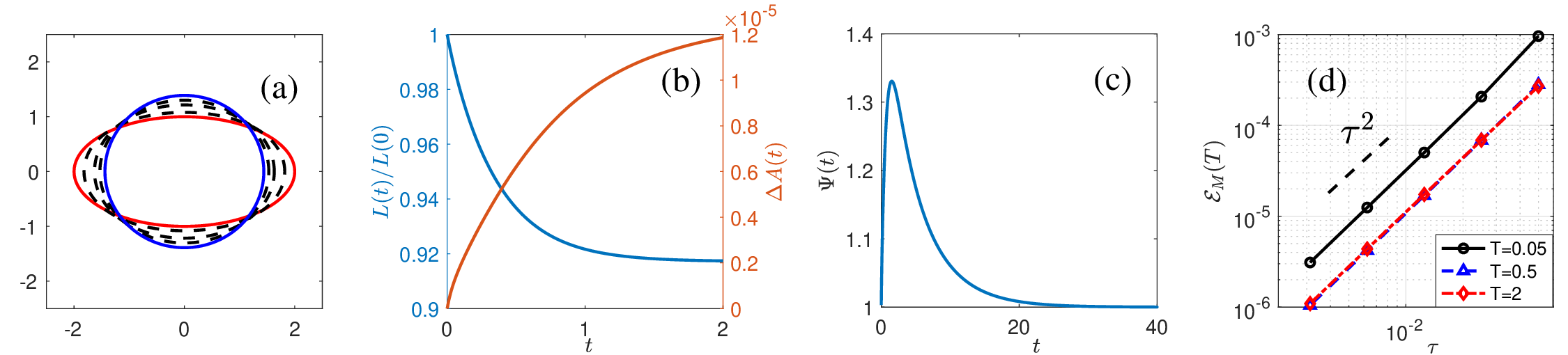}
		\caption{Evolution of an elliptic curve driven by the AP-CSF via the BGN/PC scheme: (a) snapshots of the evolution; (b) evolution of the normalized perimeter and relative area loss; (c) evolution of the mesh ratio function; (d) temporal convergence rates at $T=0.05,0.5,2$, with the spatial mesh size fixed as $N=10^4$.}
		\label{Fig:AP-CSF-BGN:PC}
	\end{figure}

The numerical results  of CSF and AP-CSF by using the BGN/PC algorithm are shown in Figures \ref{Fig:CSF} and \ref{Fig:AP-CSF-BGN:PC}, respectively.
From these Figures, it is evident that the BGN/PC schemes achieve second-order accuracy in time for both CSF and AP-CSF. Additionally, since the AP-CSF has an equilibrium state whereas the CSF does not, the mesh ratio function for AP-CSF converges to one during long-term evolution, indicating mesh equidistribution (as shown in Figure  \ref{Fig:AP-CSF-BGN:PC}(c)).

\medskip
\subsection{For anisotropic surface diffusion flow of planar curves}

 In this subsection, we apply the proposed BJL/PC scheme \eqref{aSDF:PC} to the A-SDF, which can be viewed as an $H^{-1}$-gradient flow with respect to the interfacial energy functional~\cite{Bao17b,Barrett2020,Bao2022a}
 \begin{equation}
 W=\int_{\Gamma}\gamma(\bn)\ \d s,
 \label{eqn:shape}
 \end{equation}
 where $\gamma(\bn)$ represents the surface energy density and $\bn$ is the unit normal vector of the curve. Consequently, the A-SDF can also be utilized to solve the following geometric variational problem:
 finding a closed curve that minimizes the aforementioned energy functional \eqref{eqn:shape} subject to the constraint that the area enclosed by the curve is a given constant.
 This minimization problem can alternatively be approached through a geometric method known as the Wulff construction~\cite{Wulff01, Bao17b}.
In the following, we will conduct comprehensive numerical experiments to demonstrate the superiority of our proposed BJL/PC algorithm \eqref{aSDF:PC} compared to the BJL scheme \eqref{aSDF:fisrt-order} for solving the A-SDF.
For simplicity, we will primarily consider several commonly used anisotropic surface energy densities~\cite{Bao2022a}, categorized as follows:
\begin{itemize}
\item (\textbf{Case 1}): The $k$-fold anisotropic surface energy density
\begin{align*}
	\gamma(\bn)=1+\beta\cos(k\theta),\quad \forall\, \bn=(-\sin\theta,\cos\theta)^\top\in \S^1,
\end{align*}
where $\beta \geq 0$ is a constant that controls the degree of anisotropy, and $k=2,3,4,6$ represents the order of rotational symmetry.
It is noteworthy that the surface energy exhibits strong anisotropy when $\beta>1/(k^2-1)$~\cite{Bao17}.
\smallskip
\item (\textbf{Case 2}): The Riemannian-like (BGN-like) metric surface energy density:
\begin{align*}
	\gamma(\bn)=\sum_{l=1}^M \sqrt{\bn^\top \mathbf{G}_{l}\bn},\quad \forall\, \bn\in \S^1,
\end{align*}
where $\mathbf{G}_l\in \R^{2\times 2}$, with $l=1,\ldots,M$, are  symmetric and positive definite matrices.
\smallskip
\item (\textbf{Case 3}): The regularized $l^1$-norm metric surface energy:
\[
\gamma(\bn)=\sqrt{n_1^2+\eps^2n_2^2}+\sqrt{\eps^2n_1^2+n_2^2},\quad \forall\, \bn=(n_1,n_2)^\top\in \S^1,
\]
where $0<\eps \ll 1$ is a small regularization constant.
\end{itemize}

\begin{ex}[Convergence order test] We investigate the convergence orders of the BJL scheme \eqref{aSDF:fisrt-order} and the BJL/PC scheme \eqref{aSDF:PC} for the evolution of an elliptic curve
under two types of anisotropic surface energy densities (i.e., Case 1 and Case 2).
\end{ex}

\begin{figure}[h!]
	\centering
	\includegraphics[width=15cm,height=8cm]{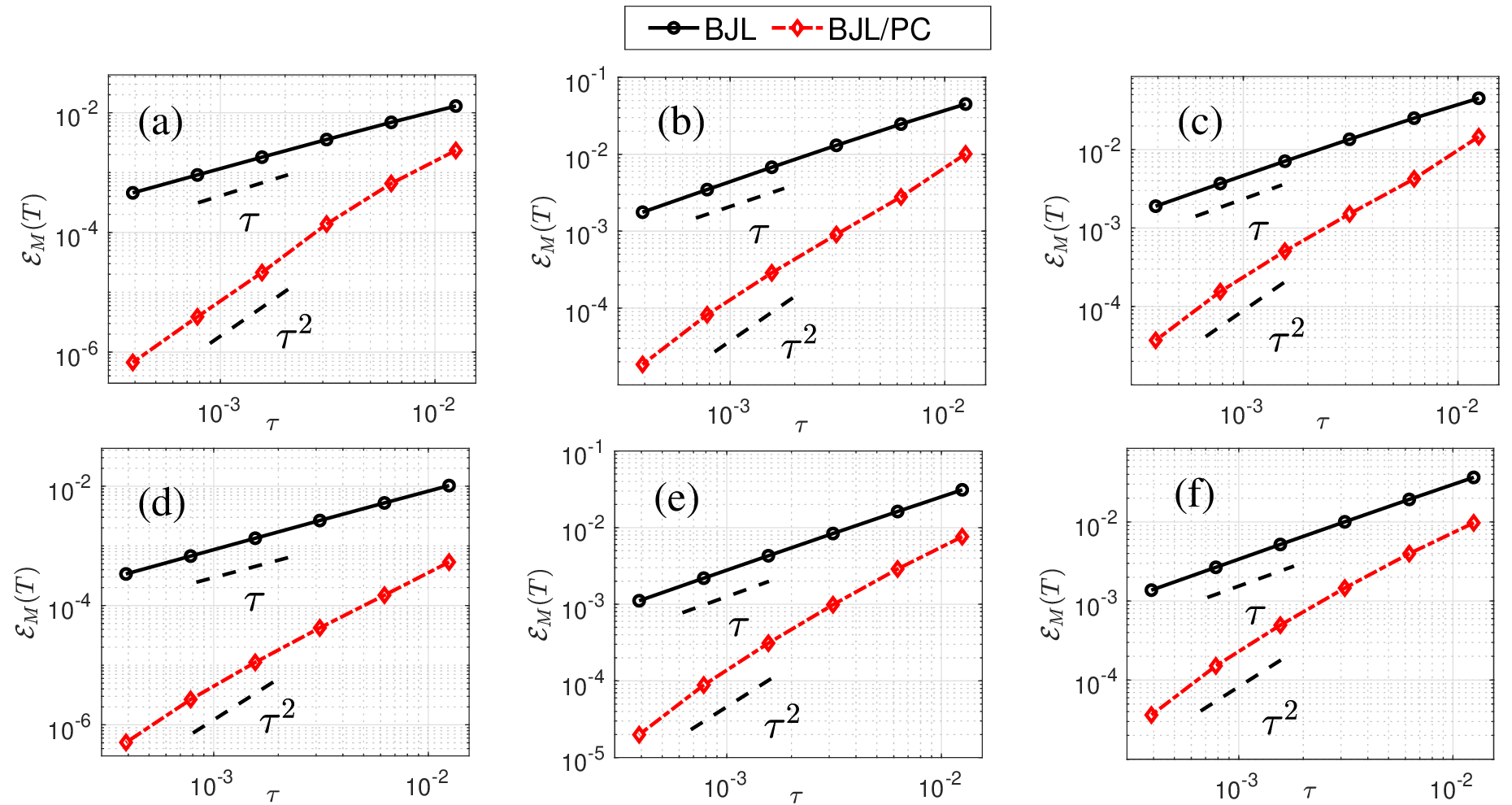}
\caption{Log-log plot of the numerical errors of the BJL scheme \eqref{aSDF:fisrt-order} and the BJL/PC scheme \eqref{aSDF:PC} for the evolution of an elliptic curve under various types of anisotropic surface energy densities:
Case 1 with $\beta=0.01$ and $k=4$ (first column); Case 1 with $\beta=0.2$ and $k=2$ (second column); Case 2 with $M=1$ and $\mathbf{G}=\mathrm{diag}(1,2)$ (last column), where the top row is for $T=0.05$ and the bottom row is for $T=0.25$.}
     \label{Fig:ani_EOC}
\end{figure}

By computing the reference solution in a similar manner to that used in Example \ref{ex:iso_EOC} and measuring numerical errors through the manifold distance,
we observe that the BJL/PC scheme \eqref{aSDF:PC} attains second-order accuracy in time (as shown in Figure \ref{Fig:ani_EOC}),
while the original BJL scheme \eqref{aSDF:fisrt-order} is only first-order accurate,
regardless of the chosen surface energy density. This comparison further underscores the superiority of our proposed BJL/PC scheme in terms of accuracy compared to the original BJL scheme.

\begin{ex}[Morphological evolution]\label{ex:Morphological evolution}
	We employ the BJL/PC scheme \eqref{aSDF:PC} to simulate the morphological evolution of an elliptic curve toward its equilibrium state under various surface energy densities.
\end{ex}

For geometric quantities, in addition to the  relative area loss $\Delta A(t)$ and the mesh ratio $\Psi(t)$, we also investigate the evolution of the normalized energy $W(t)/W(0)$, defined as
\[
    \l.\frac{W(t)}{W(0)}\r|_{t=t_m} = \frac{W^m}{W^0},\quad W^m=\sum_{j=1}^N|\mathbf{h}_j^m|\gamma(\bn_j^m),
\]
where $W^m$ represents the discrete energy of the curve $\Gamma^m$.

\begin{figure}[h!]
		\centering
		\includegraphics[width=16cm,height=10cm]{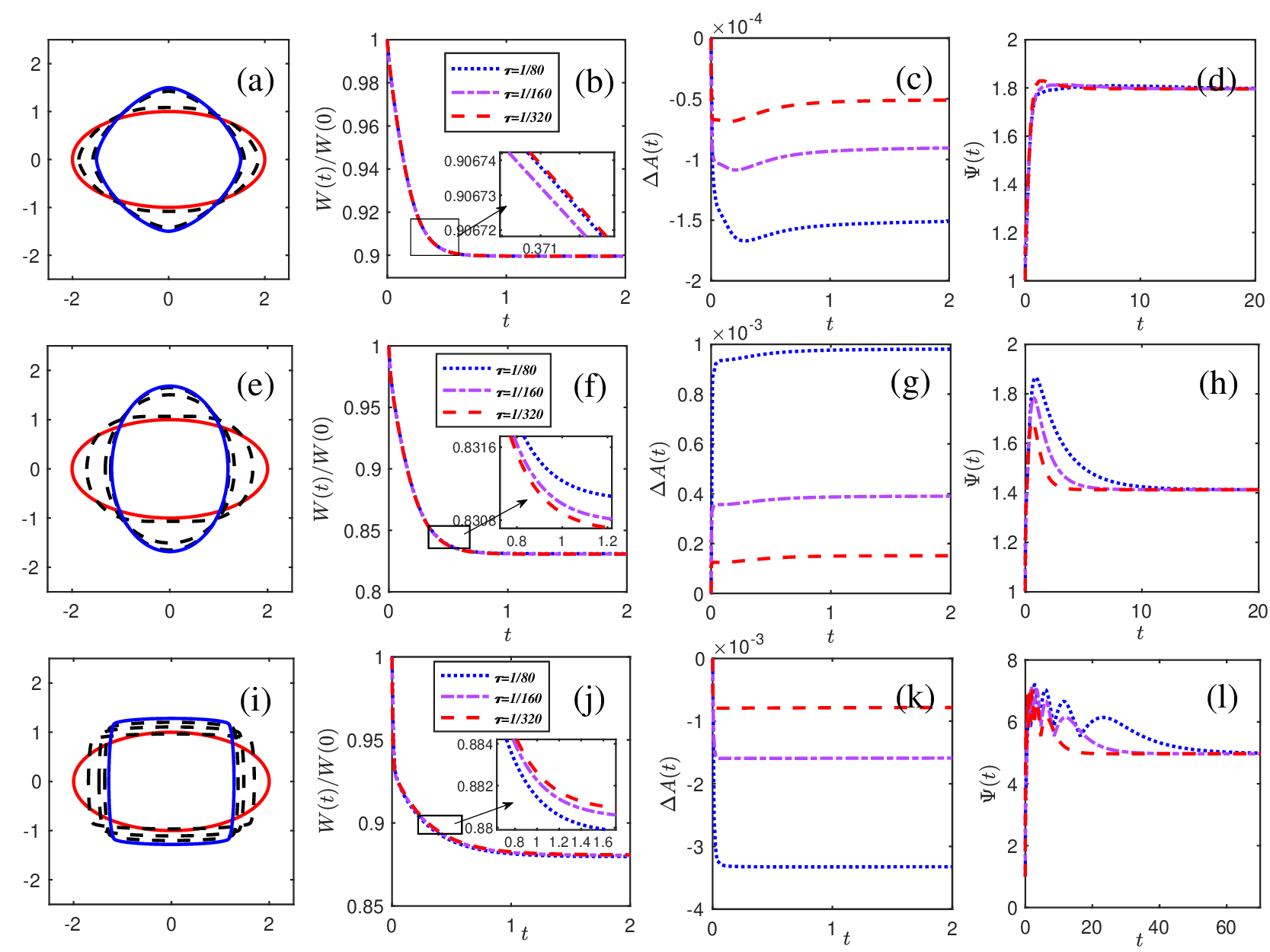}
\vspace{-1mm}		
		\caption{Evolution of curves starting from an elliptic configuration (red line) toward their equilibrium states (blue line), and the corresponding geometric quantities using the BJL/PC scheme \eqref{aSDF:PC} under three weakly anisotropic surface energy densities: Case 1 with $\beta=0.05$ and $k=4$ (top row); Case 2 with $M=1$ and $\mathbf{G}=\mathrm{diag}(1,2)$ (second row); Case 3 with $\varepsilon=0.01$ (bottom row). The parameters are set to $N=320$, and $\tau=1/80,1/160,1/320$.}
		\label{Fig:EVO_ani}
	\end{figure}

\begin{figure}[h!]
\centering
\includegraphics[width=15cm,height=7cm]{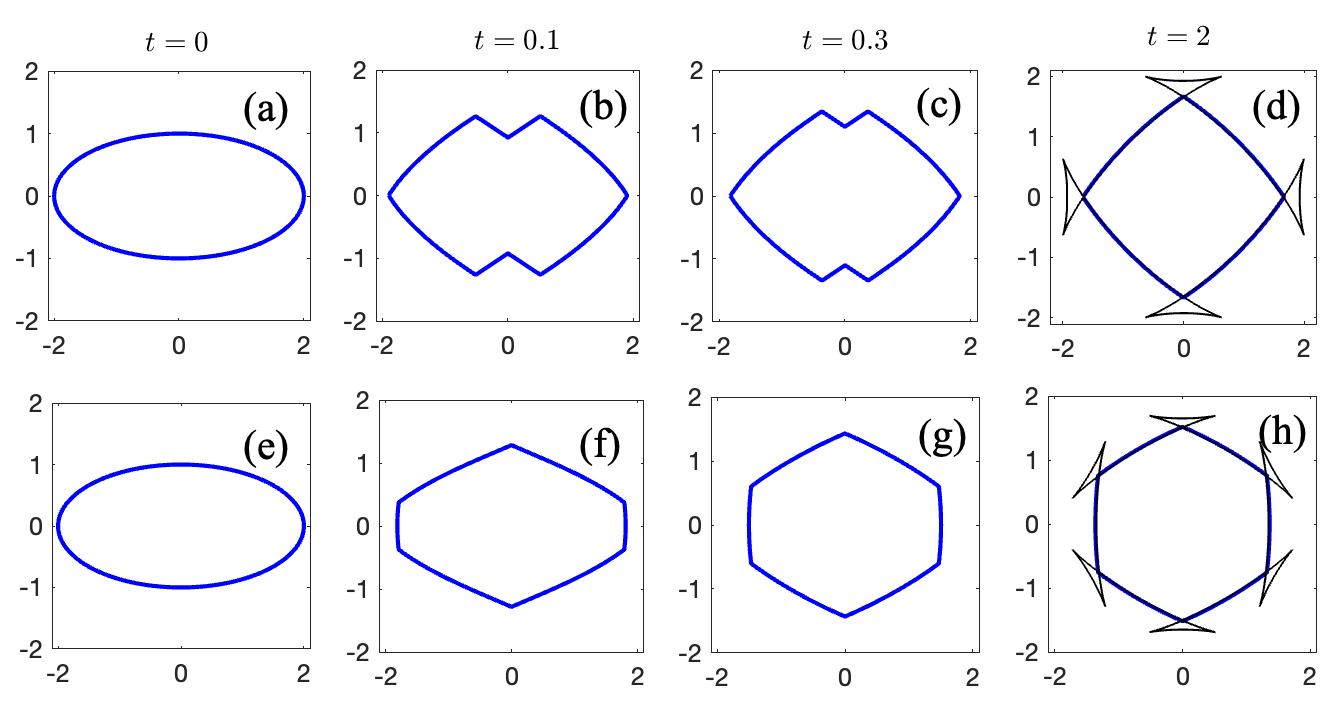}
\vspace{-5pt}
\caption{Evolution of  an initially elliptic curve driven by the A-SDF toward its equilibrium using the BJL/PC scheme \eqref{aSDF:PC} for two strongly anisotropic surface energy densities:
Case 1 with  $\beta=0.2$ and $k=4$ (top row); Case 1 with $\beta=0.1$ and $k=6$ (bottom row). The discretization parameters are set to $N=320$ and $\tau=1/320$. The black solid line, depicted in (d) and (h), represents the Wulff envelope.}
\label{Fig:EVO_s_ani}
\end{figure}

We firstly investigate the evolution of elliptic curves governed by the A-SDF with three different weakly anisotropic surface energies. As illustrated in Figure \ref{Fig:EVO_ani}, our observations reveal that the BJL/PC scheme \eqref{aSDF:PC} effectively evolves curves with three distinct surface energy types towards their respective equilibrium states. These numerically computed equilibrium states, represented as blue lines in Figures \ref{Fig:EVO_ani}(a), (e) and (i), closely align with the theoretical equilibrium obtained through the Wulff construction.
Furthermore, we explore the corresponding geometric quantities and find that the BJL/PC scheme monotonically decreases the discrete energy $\Gamma^m$. Notably, it effectively preserves the enclosed area with a numerical error margin below $0.1\%$. In anisotropic scenarios, although the BJL/PC scheme \eqref{aSDF:PC} does not attain mesh equidistribution over long-time evolution, it maintains a remarkably high level of mesh quality, as evidenced by consistently low mesh ratios throughout the entire evolution.

 Figure \ref{Fig:EVO_s_ani} presents the simulation for elliptical curves with two distinct strongly anisotropic surface energy densities. Specifically, Case 1 is characterized by either $\beta=0.2$ and $k=4$ or $\beta=0.1$ and $k=6$. In Figures \ref{Fig:EVO_s_ani}(d) and (h), the numerically computed equilibrium states (in blue) exhibit a remarkable agreement with their corresponding theoretical equilibriums derived from the Wulff construction~\cite{Wulff01, Bao17b}, which suggests that the theoretical equilibrium state is obtained by trimming the `ears' of the Wulff envelope (shown in black).

\begin{ex}[Comparison of numerical equilibrium states]
	We compare the numerical equilibrium states obtained from the BJL scheme \eqref{aSDF:fisrt-order} and BJL/PC scheme \eqref{aSDF:PC} with the theoretical equilibrium derived from the Wulff construction in both weakly and strongly anisotropic cases.
\end{ex}

\begin{figure}[h!]
\centering
\includegraphics[width=14cm,height=6.5cm]{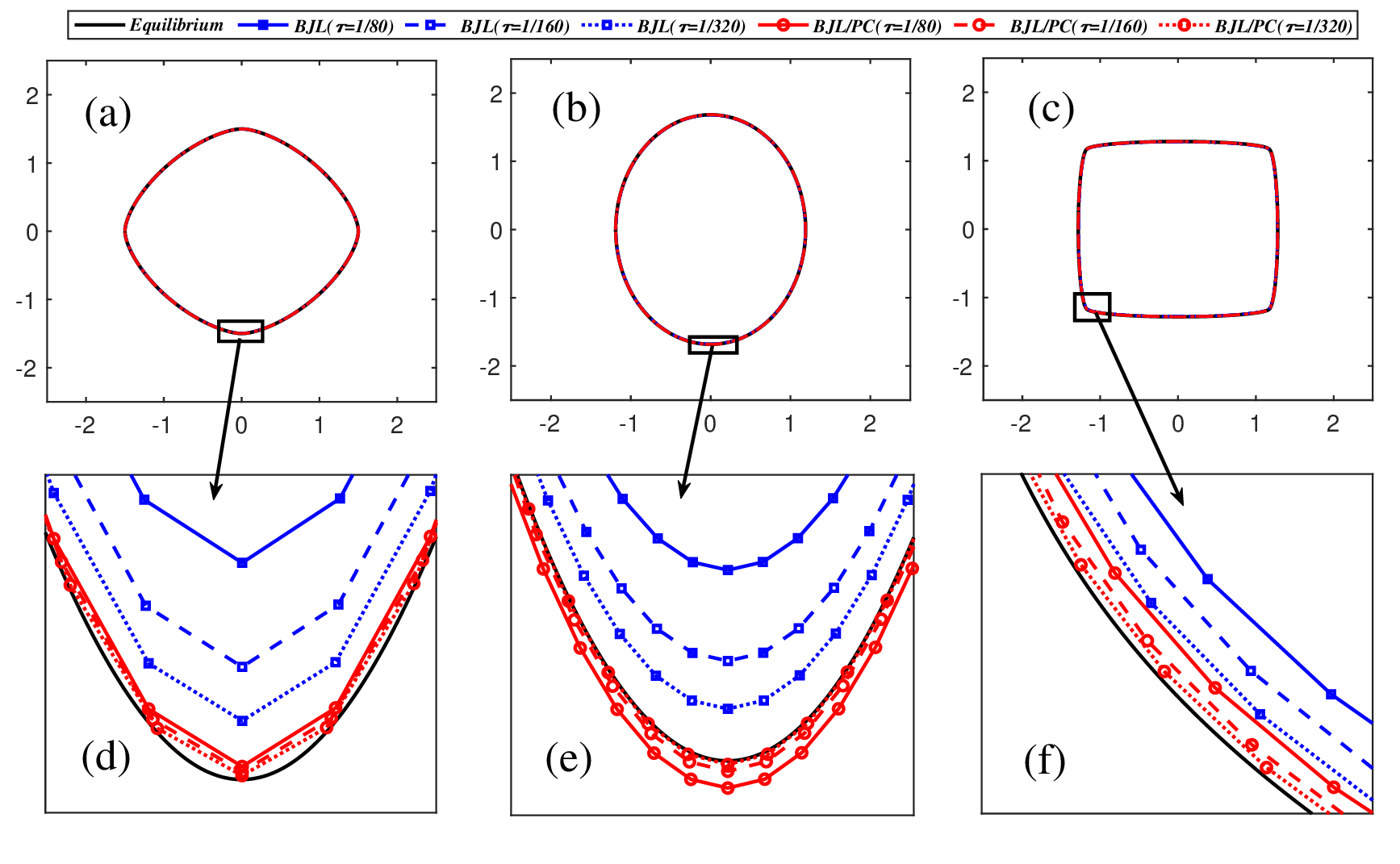}
\setlength{\abovecaptionskip}{-3pt}
\vspace{10pt}
\caption{Numerical equilibrium states obtained from the BJL scheme \eqref{aSDF:fisrt-order} (in blue) and the  BJL/PC scheme \eqref{aSDF:PC} (in red) for three weakly anisotropic surface energy densities: (a) Case 1 with $\beta=0.05$ and $k=4$; (b) Case 2 with $M=1$ and $\mathbf{G}=\mathrm{diag}(1,2)$; (c) Case 3 with $\varepsilon=0.01$. The corresponding zoom-in figures are presented in (d)-(f). The parameter are set to $N=320$ and  $\tau=1/80,1/160,1/320$. The black solid lines represent the theoretical equilibriums constructed using the Wulff  construction.}
\label{Fig:Wulff_shape}
\end{figure}

\begin{figure}[h!]
\centering
\includegraphics[width=11cm,height=7cm]{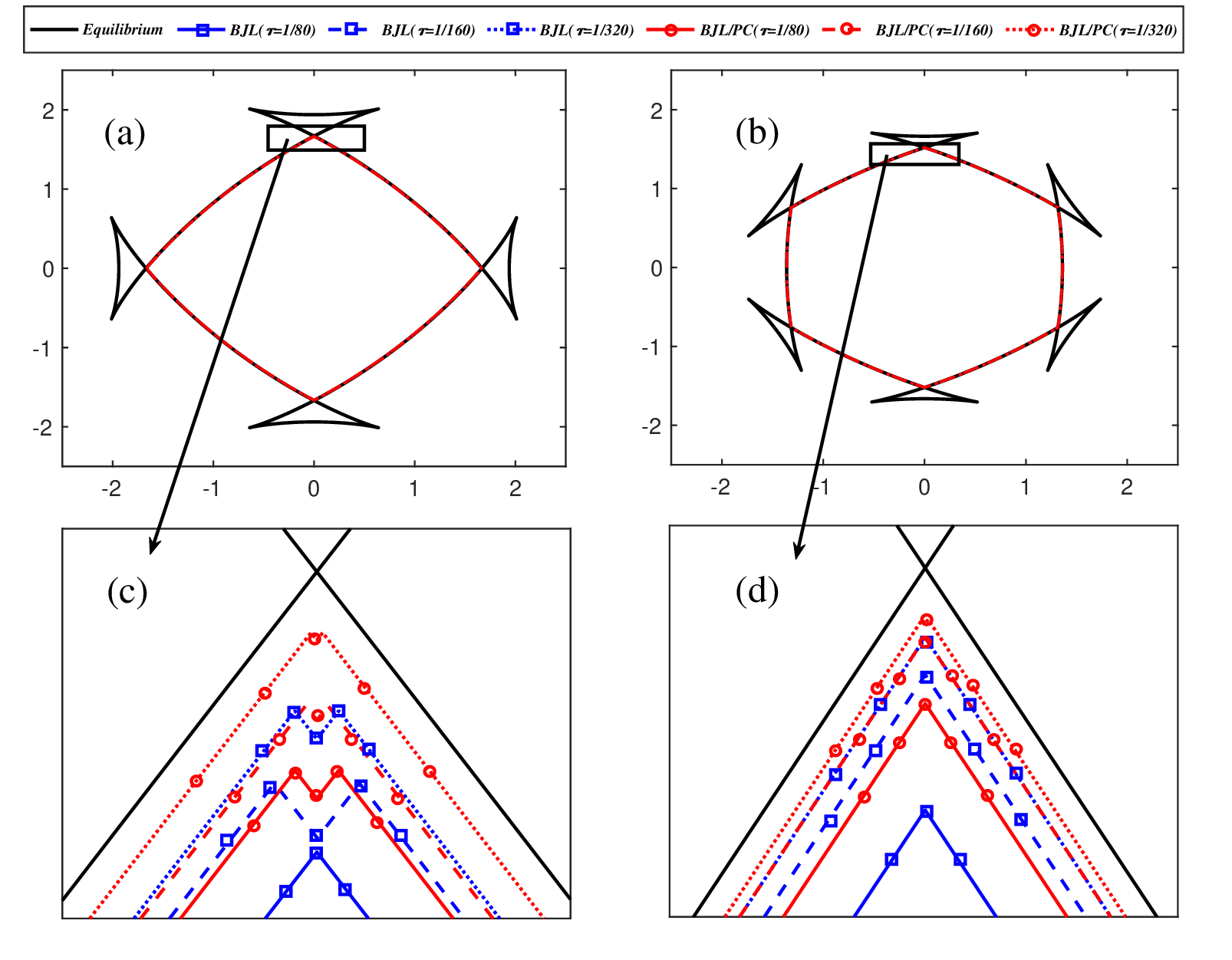}
\caption{Numerical equilibriums obtained from the BJL scheme \eqref{aSDF:fisrt-order} (in blue) and the BJL/PC scheme \eqref{aSDF:PC} (in red) for two strongly anisotropic surface energy densities: (a) Case 1 with  $\beta=0.2$ and $k=4$; (b) Case 1 with $\beta=0.1$ and $k=6$.  The corresponding zoom-in plots are displayed in (c)-(d) and the computational parameters are set to $N=320$ and $\tau=1/80,1/160,1/320$. The black solid lines represent the Wulff envelope in the strongly anisotropic case.}
\label{Fig:sharp_corner}
\end{figure}

Figures \ref{Fig:Wulff_shape} and \ref{Fig:sharp_corner} provide detailed comparisons of the equilibriums obtained using the BJL scheme \eqref{aSDF:fisrt-order} (shown in blue) and the BJL/PC scheme \eqref{aSDF:PC} (shown in red) for weakly and strongly anisotropic energy densities, respectively. The corresponding surface energy densities are specified in Figures \ref{Fig:EVO_ani} and \ref{Fig:EVO_s_ani}. For these simulations, we maintain a fixed number of mesh points at $N=320$ while
varying the time step $\tau$. The results demonstrate that the BJL/PC scheme \eqref{aSDF:PC} exhibits higher accuracy than the BJL scheme \eqref{aSDF:fisrt-order} in predicting the Wulff shape for both weakly and strongly anisotropic cases. Notably, in Figures \ref{Fig:sharp_corner}(c)-(d), the BJL/PC scheme effectively captures the sharp corners characteristic of the strongly anisotropic scenario.

The comparisons further demonstrate that the BJL/PC scheme \eqref{aSDF:PC} not only captures the overall features of the equilibrium states but also accurately reproduces finer details. This is evident from the close alignment between the numerical and theoretical equilibriums, underscoring the robustness and precision of the BJL/PC scheme in simulating the evolution of curves in the case of anisotropic surface energy densities.

\subsection{For surface diffusion flow of surfaces in $\R^3$}

 In this part, we extend the BGN/PC scheme \eqref{SDF3d:BGN/PC} to the SDF in $\R^3$. Similar to the planar case discussed earlier, we compare the numerical results against those obtained using the BGN scheme, the BGN/CNLF scheme and the BGN/BDF2 scheme. The latter two schemes are designed for surface evolution in a manner analogous to \eqref{SDF:BGN2} and \eqref{SDF:BDF2}. As noted in Remark \ref{Remark:BGN2}, mesh regularization techniques are necessary for the BGN/CNLF scheme. In practical simulations, we implement this mesh regularization for the BGN/CNLF scheme by applying the first-order BGN scheme \eqref{SDF3d:BGN1} for the trivial flow (cf. \cite[Subsection 2.4, Remark 2.3]{Jiang2023}).

\begin{ex}[Convergence order test] We compare the convergence orders of the  BGN scheme \eqref{SDF3d:BGN1}, the BGN/BDF2 scheme, and the  BGN/PC scheme \eqref{SDF3d:BGN/PC} for the evolution of a $2:1:1$ ellipsoidal surface defined by $x^2/4+y^2+z^2=1$.
\end{ex}

\vspace{-5pt}
\begin{figure}[h!]
		\centering
		\includegraphics[width=8cm,height=6cm]{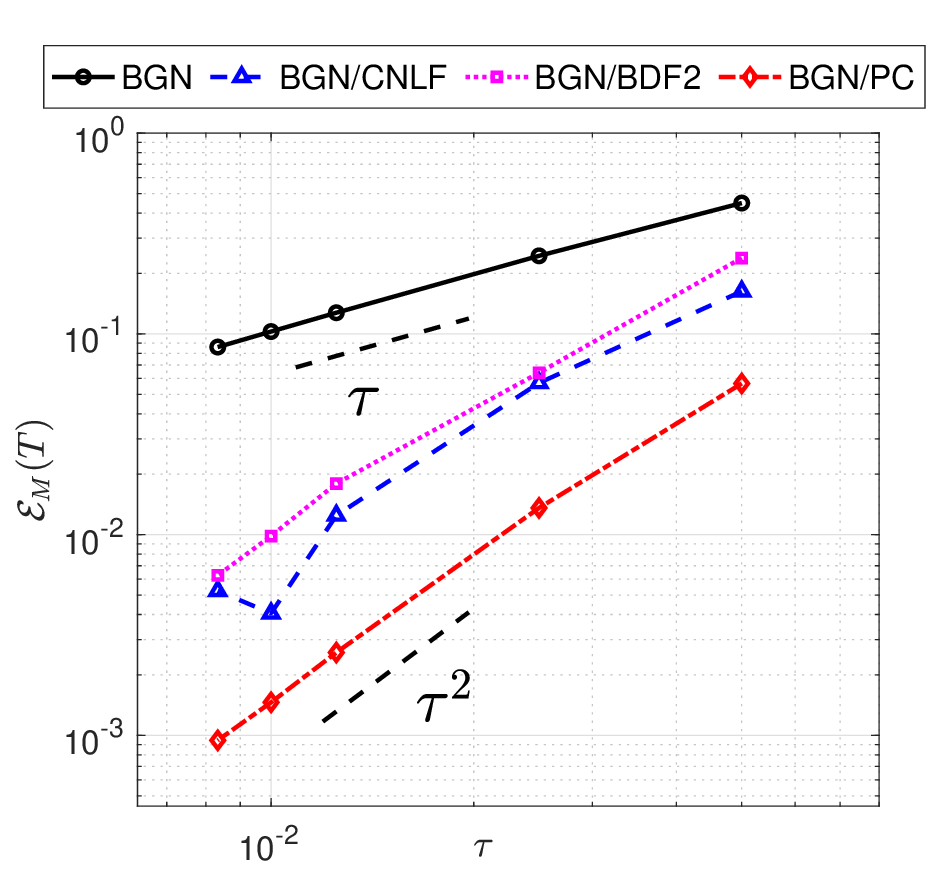}
		\vspace{-5pt}
		\caption{Numerical errors of the three schemes for solving the SDF in $\R^3$ at time $T=0.15$, with the initial surface chosen as a $2:1:1$ ellipsoid. The spatial mesh is configured as $(J,K)=(100280,50142)$.}
		\label{Fig:SDF3d_EOC}
	\end{figure}

To assess the convergence order of each method, we utilize the reference solution $\bX_{\mathrm{ref}}$ obtained from the corresponding scheme, employing an exceptionally fine mesh with parameters $(J,K)=(100280,50142)$ and a time step of $\tau=1/1500$. Here, $J$ and $K$ denote the number of triangles and vertices of the initial polyhedron, respectively. We then maintain this fine mesh to evaluate the temporal error. The numerical error and convergence order are defined analogously to the planar case. Figure \ref{Fig:SDF3d_EOC} displays a log-log plot of the manifold distance at time $T=0.15$ for the BGN scheme, the BGN/CNLF scheme, the BGN/BDF2 scheme, and the BGN/PC scheme. It is evident that the BGN and BGN/PC schemes achieve the desired first- and second-order convergence, respectively. However, there is a notable reduction in order in certain regions for the other two second-order schemes. Furthermore, the BGN/PC scheme is the most accurate among four schemes studied. These observations indicate that the BGN/PC scheme demonstrates superior robustness and accuracy compared to the other methods presented.

\begin{ex}[Evolution of morphology and geometric quantities]\label{ex:SDF3d_evo}
We apply the BGN scheme \eqref{SDF3d:BGN1}, the BGN/CNLF scheme, the BGN/BDF2 scheme, and the BGN/PC scheme \eqref{SDF3d:BGN/PC} to simulate the morphological evolution and their geometric quantities of a $2:1:1$ ellipsoidal surface and a torus defined by
	$(\sqrt{x^2+y^2}-1)^2+z^2=4/25$.
\end{ex}

\begin{figure}[h!]
\vspace{0mm}
		\centering
		\includegraphics[width=15cm,height=12cm]{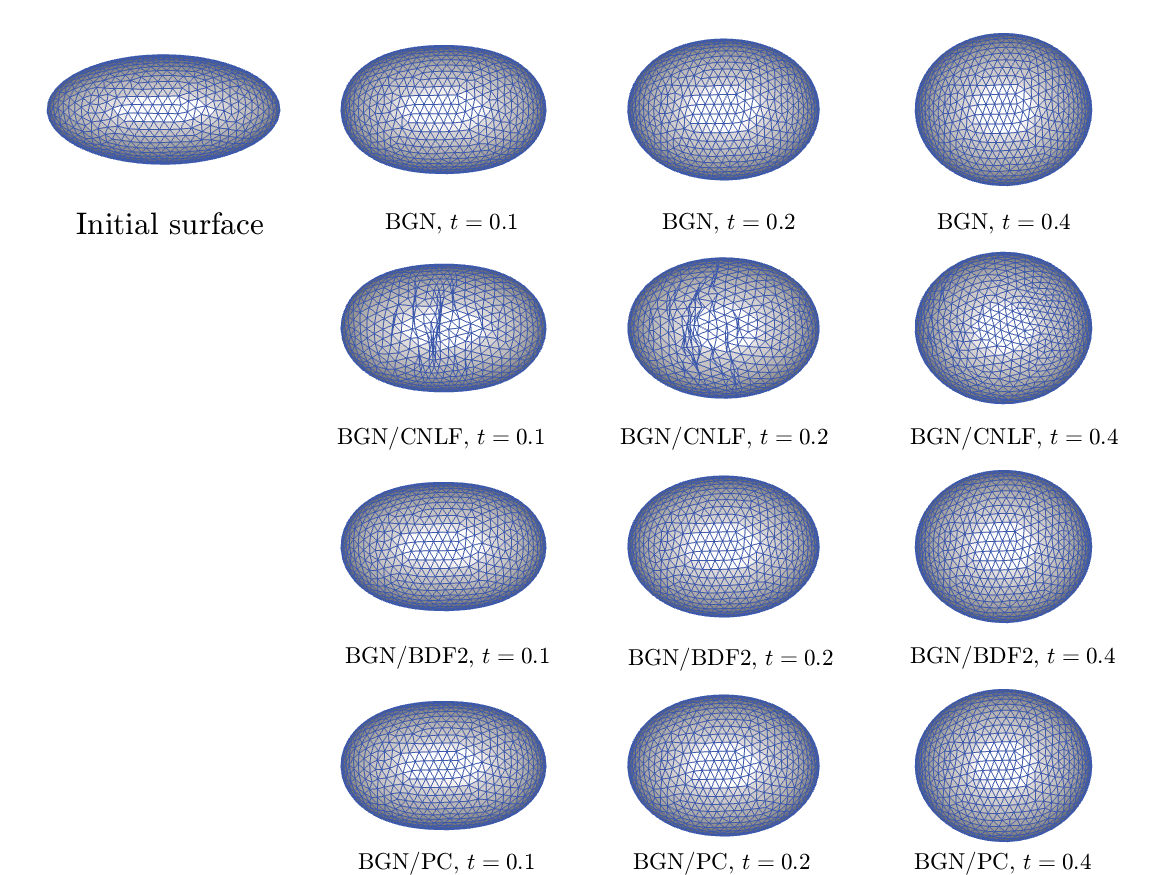}
		\caption{Evolution of a $2:1:1$ ellipsoid driven by the SDF using the BGN scheme (top row), the BGN/CNLF scheme (second row), the BGN/BDF2 scheme (third row) and the BGN/PC scheme (bottom row).
	The spatial mesh is configured as $(J,K)=(2052,1028)$ and the time step is set as $\tau=1/550$.}
		\label{Fig:SDF3d_EOC}
	\end{figure}

\begin{figure}[h!]
		\centering
	\includegraphics[width=15cm,height=4cm]{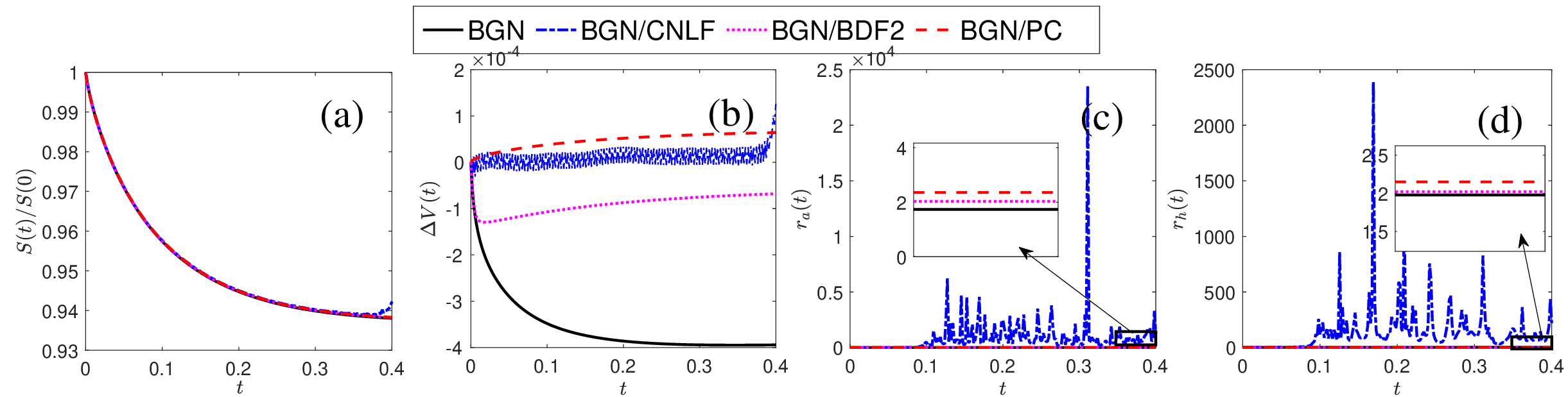}
		\caption{Evolution of the corresponding geometric quantities of the ellipsoid: (a) the normalized surface area; (b) the relative volume loss; (c)-(d) the mesh distribution functions $r_a(t)$ and $r_h(t)$.}
		\label{Fig:SDF3d_ell_geo}
	\end{figure}

\begin{figure}[h!]
\vspace{0mm}
		\centering
	\includegraphics[width=14cm,height=13cm]{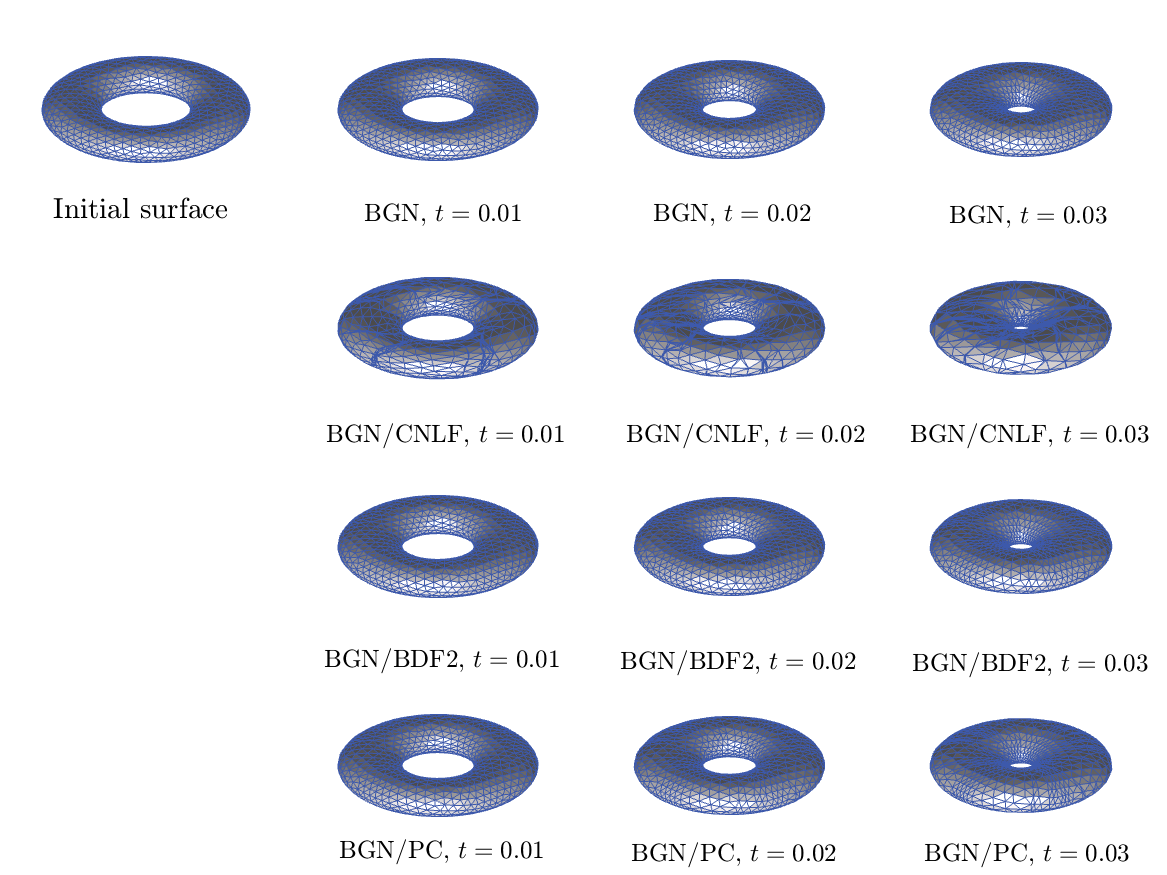}
		\caption{Evolution of a torus driven by the SDF using the BGN scheme (top row), the BGN/CNLF scheme (second row), the BGN/BDF2 scheme (third row) and the BGN/PC scheme (bottom row), with a spatial mesh configured as $(J,K)=(2000,1000)$ and a time step of $\tau=1/4000$.}
	\label{Fig:SDF3d_torus_evo}
	\end{figure}

\begin{figure}[h!]
		\centering
\includegraphics[width=15cm,height=4cm]{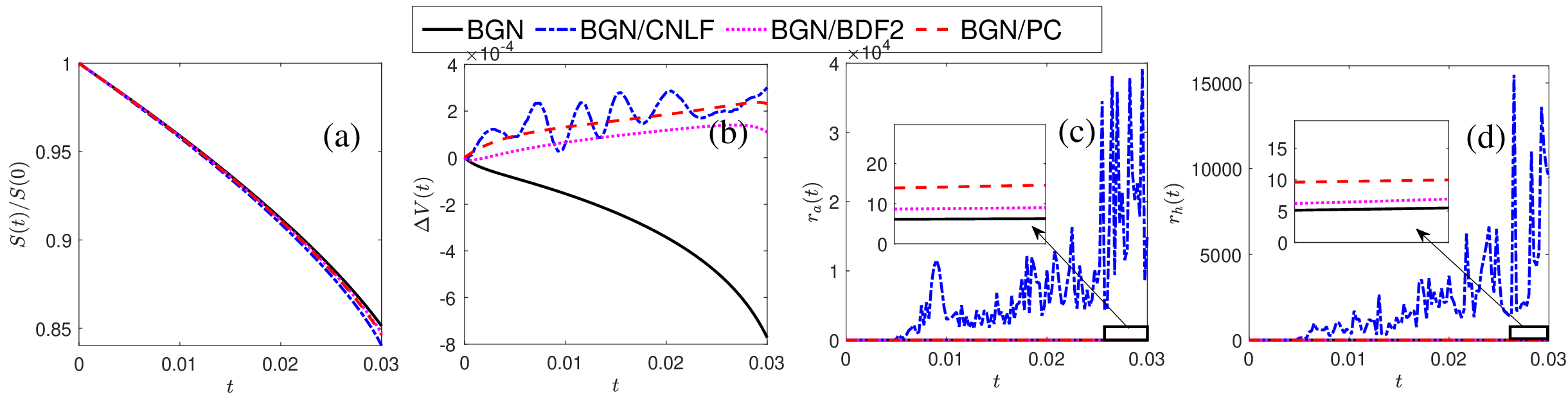}
		\caption{Evolution of the corresponding geometric quantities of the torus: (a) the normalized surface area; (b) the relative volume loss; (c)-(d) the mesh distribution functions $r_a(t)$ and $r_h(t)$.}
		\label{Fig:SDF3d_torus_geo}
	\end{figure}

\begin{figure}[h!]
		\centering
\includegraphics[width=14cm,height=4cm]{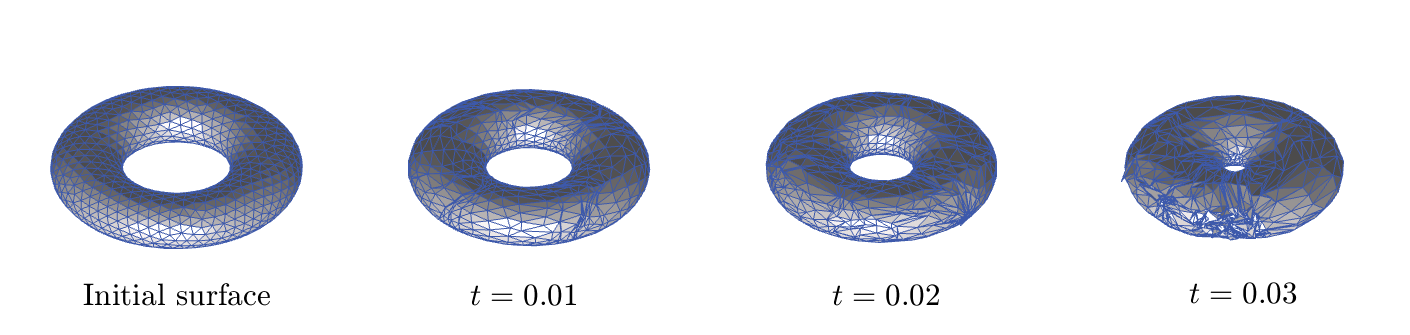}
\caption{Evolution snapshots of a torus driven by the SDF using the BGN/CNLF scheme without mesh regularization.}
	\label{toruswrong}
\end{figure}

Figures \ref{Fig:SDF3d_EOC} and \ref{Fig:SDF3d_torus_evo} illustrate the evolution of an ellipsoid and a torus using the four
distinct schemes, respectively. As clearly shown in these plots, the two second-order schemes, the BGN/BDF2 scheme and the BGN/PC scheme, exhibit highly stable numerical evolutions without requiring any mesh regularization,
similar to the first-order BGN scheme. In contrast, the BGN/CNLF scheme encounters mesh distortion and necessitates mesh regularization, as illustrated in the second row of Figure \ref{Fig:SDF3d_torus_evo} and Figure \ref{toruswrong}.

Additionally, we investigate the evolution of the following geometric quantities, as shown in Figures \ref{Fig:SDF3d_ell_geo} and \ref{Fig:SDF3d_torus_geo}:
(1) the relative volume loss $\Delta V(t)$; and (2) the normalized surface area $S(t)/S(0)$, which are defined as follows for $m\ge 0$:
\[
\Delta V(t)|_{t=t_m}=\frac{V^m-V^0}{V^0},\quad \l.\frac{S(t)}{S(0)}\r|_{t=t_m}=\frac{S^m}{S^0},
\]
where $V^m$ and $S^m$ represent the volume enclosed by the polyhedron $\Gamma^m$ and the surface area of $\Gamma^m$, respectively. Furthermore, two mesh distribution functions $r_h(t)$ and $r_a(t)$, are defined to characterize mesh quality~\cite{Jiang2024}:
\[
	\l.r_h(t)\r|_{t=t_m}
	:=\frac{\max_{j}\max \{\|\mathbf{q}^m_{j_1}-\mathbf{q}^m_{j_2}\|,\|\mathbf{q}^m_{j_2}-\mathbf{q}^m_{j_3}\|,\|\mathbf{q}^m_{j_3}-\mathbf{q}^m_{j_1}\| \}}{\min_{j}\min \{\|\mathbf{q}^m_{j_1}-\mathbf{q}^m_{j_2}\|,\|\mathbf{q}^m_{j_2}-\mathbf{q}^m_{j_3}\|,\|\mathbf{q}^m_{j_3}-\mathbf{q}^m_{j_1}\|\}},\quad
\l.r_a(t)\r|_{t=t_m}
	:=\frac{\max_j|\sigma_j^m|}{\min_j|\sigma_j^m|}.
\]

Figures \ref{Fig:SDF3d_ell_geo} and \ref{Fig:SDF3d_torus_geo} (a)-(b) demonstrate that both the BGN/PC scheme and the BGN/BDF2 scheme are more effective in preserving geometric properties, specifically with regard to area reduction and volume conservation. Notably, the volume loss observed in these two second-order schemes is significantly lower than that in the first-order BGN scheme. In the case of the BGN/CNLF scheme, oscillations are evident in both volume loss and mesh ratio functions. Furthermore, Figures \ref{Fig:SDF3d_ell_geo} and \ref{Fig:SDF3d_torus_geo} (c)-(d) clearly illustrate significant mesh distortion for both initial surfaces when employing the BGN/CNLF method. Although mesh regularization is implemented at each time step, the mesh quality remains unsatisfactory, in contrast to the results observed in the evolution of curves \cite{Jiang2023}. Nonetheless, It is important to emphasize that mesh regularization is essential in the BGN/CNLF scheme; without it, the mesh becomes severely distorted, resulting in inaccurate numerical solutions, as illustrated in Figure \ref{toruswrong}. In contrast, the mesh distribution functions $r_h(t)$ and $r_a(t)$ of the other three schemes maintain favorable bounds throughout the evolution.

%

\section{Conclusions}

We proposed a novel parametric finite element method based on the BGN framework, which is second-order accurate in time and tailored specifically for solving surface diffusion flows. This method employs a predictor-corrector approach for time discretization, enhancing the standard first-order BGN scheme to achieve second-order accuracy. Unlike previous second-order methods, our proposed scheme maintains excellent mesh quality without requiring mesh regularization. Furthermore, our algorithm can be easily extended to address various geometric flows, including curve-shortening flow, area-preserving curve-shortening flow, anisotropic surface diffusion flow, and surface diffusion flow in $\R^3$. To validate our approach, we presented extensive numerical results demonstrating that the predictor-corrector BGN-based schemes achieve quadratic convergence in time and outperform the existing schemes.

Additionally, we emphasize that our predictor-corrector strategy can be easily adapted to solve Willmore flows \cite{Barrett2008b} or anisotropic surface diffusion flows in $\R^3$~\cite{Bao2024}.
However, the design of an unconditionally energy-stable, temporally high-order scheme for isotropic or anisotropic surface diffusion flows, as well as other geometric flows, presents a compelling and significant challenge. In future work, we intend to further investigate high-order, structure-preserving numerical schemes that effectively maintain the geometric properties of various geometric flows \cite{Garcke2024, Jiang2022}.

\section*{Acknowledgement}
Wei Jiang and Lian Zhang are supported by the National Natural Science Foundation of China (No. 12271414) and the Guangdong Basic and Applied Basic Research Foundation (No. 2024A1515012505). Chunmei Su and Ganghui Zhang are supported by National Key R\&D Program of China (No. 2023YFA1008902) and  National Natural Science Foundation of China (No. 12201342).
The numerical calculations in this paper have been partially done on the Center of High Performance Computing, Tsinghua University and the supercomputing system in the Supercomputing Center of Wuhan University.

\bigskip


\end{document}